\documentclass[letterpaper, paper,10.5pt]{AAS}		% for final proceedings (20-page limit)

\usepackage{bm}
\usepackage{amsmath}
\usepackage{amssymb}
\usepackage{overcite}
\usepackage{footnpag}			      	% make footnote symbols restart on each page
\usepackage{acronym}
\usepackage{graphicx}
\usepackage{caption}
\usepackage{subcaption}
\usepackage{siunitx}
\usepackage{tabularx}
\usepackage{booktabs}%
\usepackage{hyperref}
\usepackage{verbatim}
\usepackage{algorithm}
\usepackage{algpseudocode}
\usepackage{mathtools}
\usepackage[thinc]{esdiff}
\usepackage{capt-of}
\usepackage{longtable,tabularx}
\usepackage[version=4]{mhchem}
\usepackage{siunitx}
\usepackage{caption}
\usepackage{subcaption}
\usepackage{arydshln}

\usepackage{enumitem}
\usepackage{tikz}
\usetikzlibrary{shapes, arrows}
 \usetikzlibrary{calc} 
 \usetikzlibrary{fit}

\newcommand{\mat}{\begin{bmatrix}}
\newcommand{\matf}{\end{bmatrix}}
\newcounter{equationset}
\PaperNumber{ (Preprint) 23-281}

\begin{document}

\title{Convex Optimization-Based Model Predictive Control for the Guidance of Active Debris Removal Transfers}

\author{M. Wijayatunga\thanks{PhD Candidate, Te P\=unaha \=Atea Space Institute, University of Auckland, Auckland, New Zealand} , R. Armellin \thanks{Professor, Te P\=unaha \=Atea   Space Institute, University of Auckland, Auckland, New Zealand} ,
H. Holt \thanks{Research Fellow, Te P\=unaha \=Atea Space Institute, University of Auckland, Auckland, New Zealand}, C. Bombardelli \thanks{Professor, Universidad Politécnica de Madrid, Spain}, L. Pirovano \thanks{Research Fellow, Te P\=unaha \=Atea Space Institute, University of Auckland, Auckland, New Zealand}} %% ordering ?

\maketitle		
\begin{abstract}
Active debris removal (ADR) missions have garnered significant interest as means of mitigating collision risks in space. This work proposes a convex optimization-based model predictive control (MPC) approach to provide guidance for such missions. While convex optimization can obtain optimal solutions in polynomial time, it relies on the successive convexification of nonconvex dynamics, leading to inaccuracies. Here, the need for successive convexification is eliminated by using near-linear Generalized Equinoctial Orbital Elements (GEqOE) and by updating the reference trajectory through a new split-Edelbaum approach. The solution accuracy is then measured relative to a high-fidelity dynamics model, showing that the MPC-convex method can generate accurate solutions without iterations.
\end{abstract}

\section{Introduction}
The space environment in the low Earth orbit (LEO) is becoming increasingly congested with space debris. As a result, the average rate of debris collisions has increased to four or five objects per year \cite{Maestrini2021}. As satellites become increasingly essential to daily life, more are added to expand space-enabled services. However, additional launches increase the risk of collision for all satellites as they further saturate space with objects, endangering critical space infrastructure. \par

Active debris removal (ADR) is the process of removing derelict objects from space, thus minimizing the build-up of unnecessary objects and lowering the probability of on-orbit collisions that can fuel a \lq\lq collision cascade \rq\rq  \cite{bonnal2013active,liou2010parametric}. ADR has become significant over the past two decades, leading to numerous studies and implementations of potential debris removal missions and technologies. The ELSA-d mission designed by Astroscale was launched in March 2021 and has successfully tested {both} rendezvous algorithms required for ADR and a magnetic capture mechanism to remove objects carrying a dedicated docking plate at the end of their missions  \footnote{\url{https://astroscale.com/elsa-d-mission-update/}}. The RemoveDebris mission by the University of Surrey is another project that demonstrated various debris removal methods, including harpoon and net capture \cite{Forshaw2020TheAS}. The CleanSpace-1 mission by the European Space Agency (ESA) aims to deorbit a~112~kg upper stage of {a} Vega rocket \footnote{\url{https://www.esa.int/Space_Safety/ClearSpace-1}}.\par

While individual removals are essential milestones towards ADR implementation, large-scale missions that target multiple objects might be necessary to compete with the current debris growth rate \cite{liou2010parametric, white2014}. Such missions are expected to rely mainly on autonomy to reduce costs while maintaining reliability. However, enabling technologies for such missions are still under development \cite{surveydebris}.  \par 
In our previous work \cite{ADRMW}, a preliminary trajectory optimization tool (PMDT) was developed to obtain time and fuel optimal trajectories of multi-ADR missions while considering the effect of $J_2$, drag, eclipses, and duty cycle. This PMDT obtained near-optimal trajectories for removing up to three debris with computational times under a minute. The PMDT utilized the classical Edelbaum's method \cite{2} to first calculate the optimal time of flight and fuel expenditures of a single transfer, then made additions to Edelbaum's method to consider the effects of atmospheric drag, engine duty cycle, and solar eclipses. Lastly, a RAAN matching algorithm that does not use fuel to perform RAAN adjustments was implemented to make the transfers cheaper. This was done by introducing an intermediate phasing orbit where the spacecraft can utilize the effect of $J_2$ perturbations to reach a desired RAAN. The PMDT allowed the optimization of fuel or time consumption of multi-debris transfers through the selection of the phasing orbits and the launch time of the mission. \par 
Also in [\citenum{ADRMW}], PMDT results were used as a reference for providing guidance for complex multi-ADR mission profiles via common guidance laws such as the $\Delta v$ law\cite{slim} and Q-law \cite{Petropoulos2005}, which showed that the PMDT reference could be tracked relatively well by these guidance laws. However, these guidance laws required more fuel than estimated by the PMDT and were only able to achieve an accuracy of $\sim$10 km in the semi-major axis and $\sim$0.1 deg in inclination. This was due to the limited accuracy inherent to the classical guidance laws, which rely heavily on approximates and simplifications. The maximum rate of change approximation of the Gauss variational equations used in the Q-law \cite{Petropoulos2005} and the simplified $\Delta v$ estimation formulae used in the $\Delta v$-law \cite{slim} provide examples.  Furthermore, classical guidance laws do not have a fail-safe - they have no means of recovery when diverged significantly from the target.
\par

In contrast, Model Predictive Control (MPC) can provide significantly more accurate control, as it can account for perturbations in real time\cite{MPC1, MPC2}. Hence, as a next step, this work focuses on providing autonomous guidance for such missions using a novel MPC method. A split Edelbaum reference is first calculated using the method discussed in our previous work \cite{ADRMW}, which is loosely tracked by a convex-based optimization in predefined uniform segments of time. If a significant deviation occurs between the real trajectory and the reference, the reference is recomputed from the current position to the target (again using the split Edelbaum method). The use of MPC in the context of spacecraft rendezvous has been recently explored by L. Ravikumar \cite{RAVIKUMAR2020518}, C. Bashnick \cite{Bashnick} and R. Vazquez \cite{VAZQUEZ2015251}, who have all identified solving optimization problems at each of the control intervals as a computationally complex and time-consuming process. Convex optimization is appealing here as it can obtain optimal solutions in polynomial time. The use of convex optimization within MPC in the context of spacecraft trajectory optimization has also been explored frequently in the literature \cite{gao2018convex,sun2018convex,wang2017convex, hu2020convex}.  However, convex optimization relies on the convexification of nonconvex dynamics and constraints, leading to the need for successive convexification and potentially inaccurate solutions when the convex representation differs strongly from the reality \cite{okelly2016non,aglietti1998non}. In this work, the accuracy, robustness, and real-time capabilities of convex-MPC guidance are enhanced through; (a) Improving the convex optimization to avoid successive iterations and inaccurate outcomes using a near-linear representation of dynamics (Generalised Equinoctial Coordinates \cite{claudio}). (b) Relaxing intermediate tracking constraints to avoid unnecessary consumption of fuel in the initial stages of the mission. (c) Recomputing the reference upon significant deviation to avoid infeasibility. These will be addressed in detail in the following sections.\par
The remainder of this paper is organized as follows. The subsequent section outlines the overview of the ADR mission that the guidance is designed for. Then, the proposed MPC guidance process is described in depth. The following section gives results for a fuel-optimal tour to remove a discarded rocket body. The spacecraft is propagated in high and low-fidelity dynamics and with different thrust error settings to fully test the MPC guidance. Finally, conclusions drawn from the study are discussed.

\section{Mission Overview}
\begin{figure}[hbt!]
    \centering
    \includegraphics[width = 0.8\textwidth]{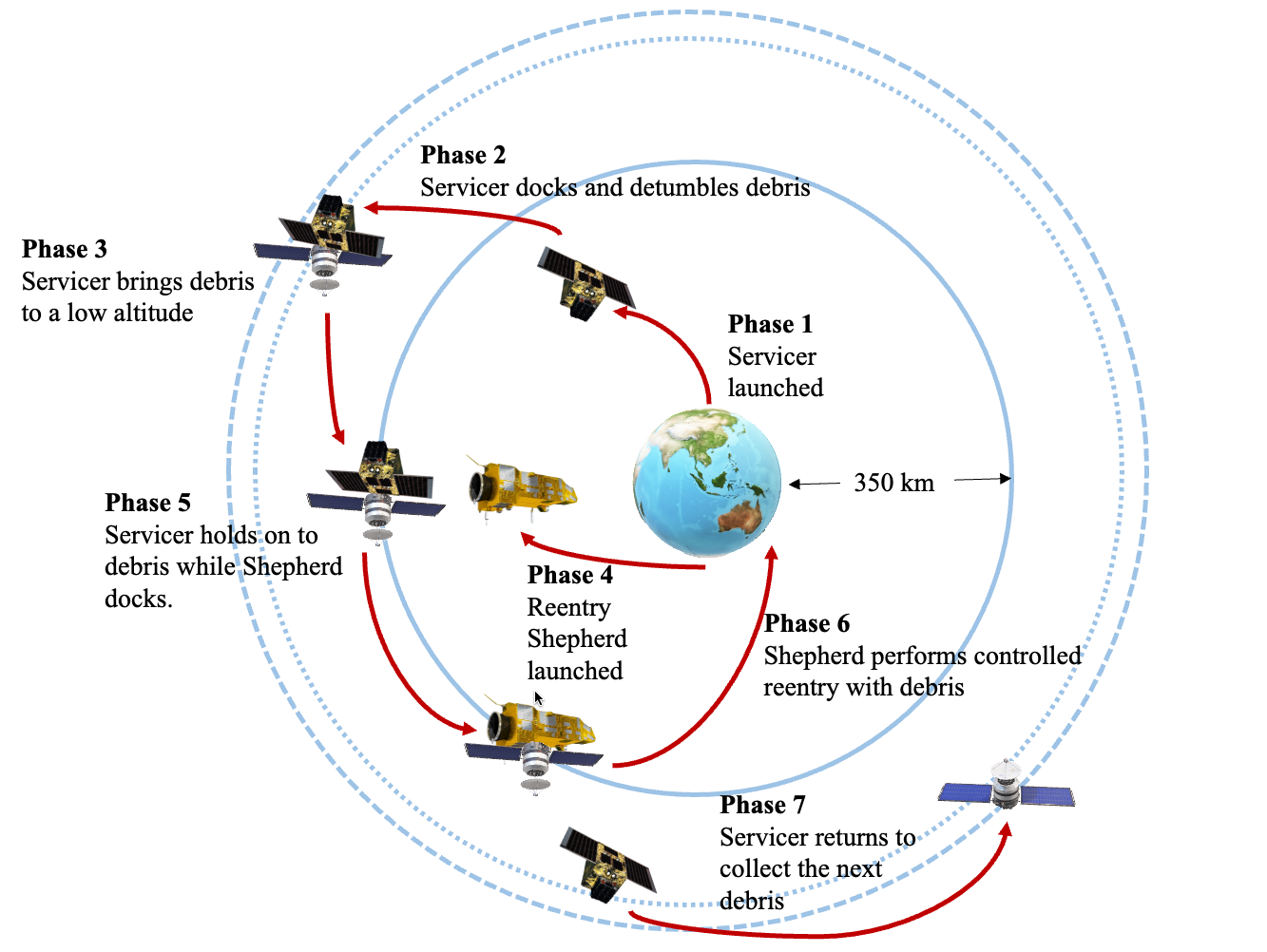}
    \caption{Mission architecture of the multi-ADR mission}
    \label{ADR}
\end{figure}

The proposed mission architecture is shown in Figure \ref{ADR}.  In this mission, two spacecraft are involved in the debris removal process. A Servicer is used to approach and rendezvous with the debris. When the rendezvous is achieved, the Servicer brings the object to a low altitude orbit ($350$~km in this study). The debris is then handed over to a Reentry Shepherd, which docks with it and performs a controlled reentry on its behalf. Controlled reentry reduces the casualty risk posed by removing the debris, which is desirable because the ADR targets are, by definition, large and thus contain components likely to survive the reentry. The Servicer {can} be reused for several debris removals, while each Reentry Shepherd can only be used once as it burns while deorbiting the debris. {A handover altitude of 350~km was selected to reduce the $\Delta v$ required of the Servicer by minimizing the orbital transfers it requires to perform while ensuring the technical feasibility of the Shepherd and satisfying safety constraints posed by the altitude of the International Space Station (ISS).}
\par 
% \par

The proposed mission can perform multi-ADR services significantly cheaper than using a single spacecraft which would perform all the mission phases and then burn in the atmosphere after removing a~single debris \cite{astroscale}. This is because the development and operation costs of a space mission are proportional to the system's dry mass~\cite{Jones2015}. Thus, although launching one Servicer and $n$ Shepherds to remove $n$ debris requires one more launch than using $n$ individual spacecraft, it is expected to be cheaper because of the lower overall mass that must be launched.
Note that this mission was previously discussed in [\citenum{ADRMW}], which details the preliminary mission design process that was conducted as a collaborative study between the University of Auckland, Astroscale and Rocketlab.  \par 
In this paper, providing guidance for the servicer spacecraft going from the initial 350 km orbit to the debris (called up leg) and coming back to the 350 km with the debris (called down leg) is studied. However, the developed guidance algorithm can be applied to any part of the mission.

\section{MPC Guidance Methodology}
The framework of the guidance methodology proposed is summarized in Figure \ref{1}. A reference to go from the initial position to a target semi-major axis ($a_f$), inclination ($i_f$) and RAAN ($\Omega_f$) is first determined using the PMDT. Then, convex-based tracking is utilized to remain as close to this trajectory as possible. This tracking is done adaptively, as the required tracking accuracy is determined by minimizing an estimate of the $\Delta v$ required to remain close to the reference. If the spacecraft deviates significantly from the reference, a new reference from the current position to the target is calculated, marking the final step of the MPC architecture.

%\comment{Need to add more time-related info here to the tikz}
\tikzstyle{connector} = [draw, -latex']
\tikzstyle{terminator} = [rectangle, draw, text centered, rounded corners, minimum height=2em]
\tikzstyle{process} = [rectangle, draw, text centered, minimum height=2em]
\tikzstyle{decision} = [diamond, draw, text centered, minimum height=2em]
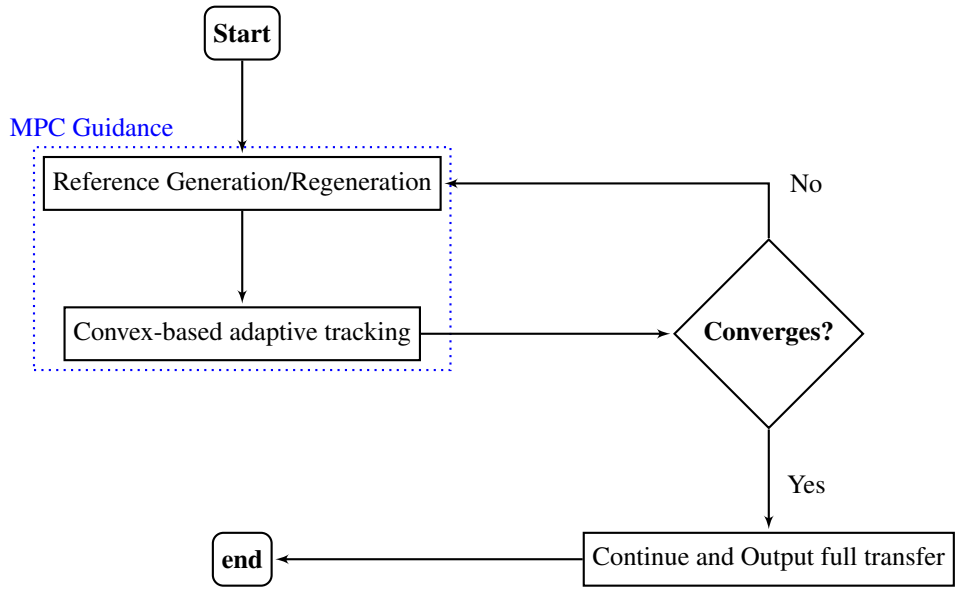
\begin{figure}[hbt!]
    \centering

\begin{tikzpicture}[->,>=stealth,shorten >=1pt,auto,node distance=2.5cm,
                    thick,main node/.style={rectangle,draw,font=\sffamily\Large\bfseries}]
      
  % Step nodes
  \node [terminator, fill=white!20] at (0,0) (start) {\textbf{Start}};
  \node[draw=none, blue] at (-2, -1.25) (MPCG) {MPC Guidance};
  \node [process, fill=white!20] at (0,-2) (step1) {Reference Generation/Regeneration};
%  \node [process, fill=white!20] at (0,-4) (step2) {Forward propagate the reference PMDT.};
    \node [process, fill=white!20] at (0,-4) (step3) {Convex-based adaptive tracking};
    \node [decision, fill=white!20] at (7,-4) (decision) { \textbf{Converges?}};
      \node [process, fill=white!20] at (7,-7) (step4) {Continue and Output full transfer};
        \node [terminator, fill=white!20] at (0,-7) (end) {\textbf{end}};
\node[draw=none] at (7.5, -2) (no) {No};
\node[draw=none] at (7.5, -6) (yes) {Yes};

 \node [fit=(step1) (step3),draw,dotted,blue] {};

\path [connector] (start) -- (step1);
%\path [connector] (step1) -- (step2);
\path [connector] (step1) -- (step3);
\path [connector] (step3) -- (decision);
\path [connector] (decision) -- (step4);
\path [connector] (decision) |- (step1);
\path [connector] (step4) -- (end);

\end{tikzpicture}
\caption{High-level overview of the Guidance procedure}
\label{1}
\end{figure}

\subsection{Reference Generation}

\begin{figure}[hbt!]
  \centering
  \begin{subfigure}[b]{0.49\textwidth}
  \centering
    \includegraphics[width = \textwidth]{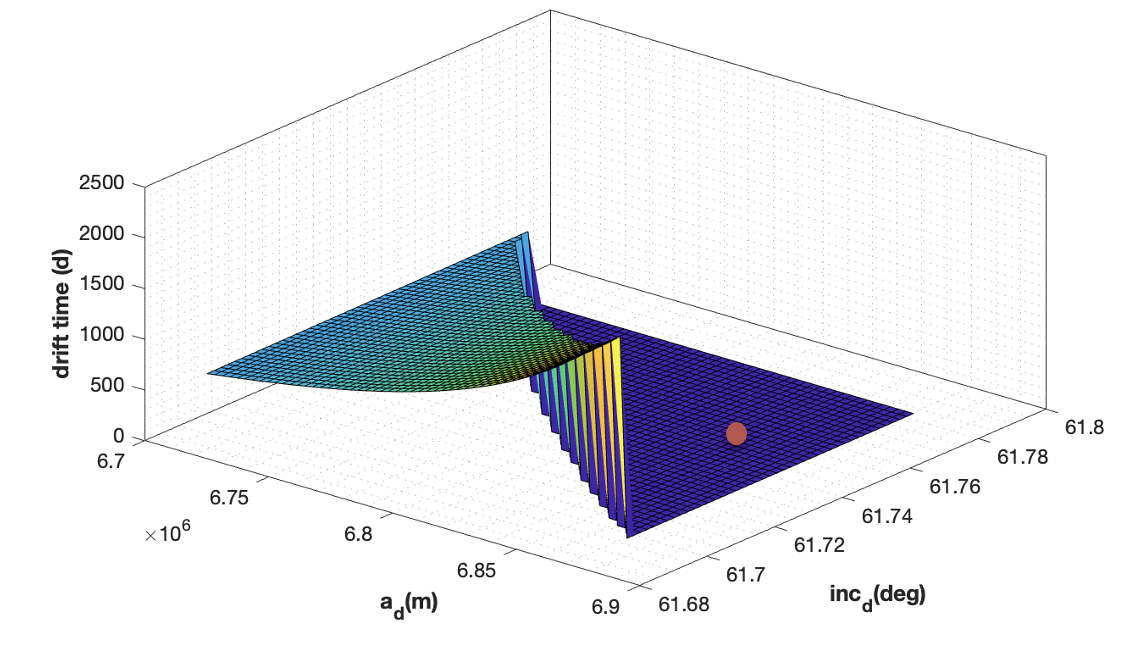}
    \caption{The discontinuous search space encountered during RAAN matching }
    \label{issue1}
  \end{subfigure}
  \hfill
  \begin{subfigure}[b]{0.49\textwidth}
    \centering
    \includegraphics[width = \textwidth]{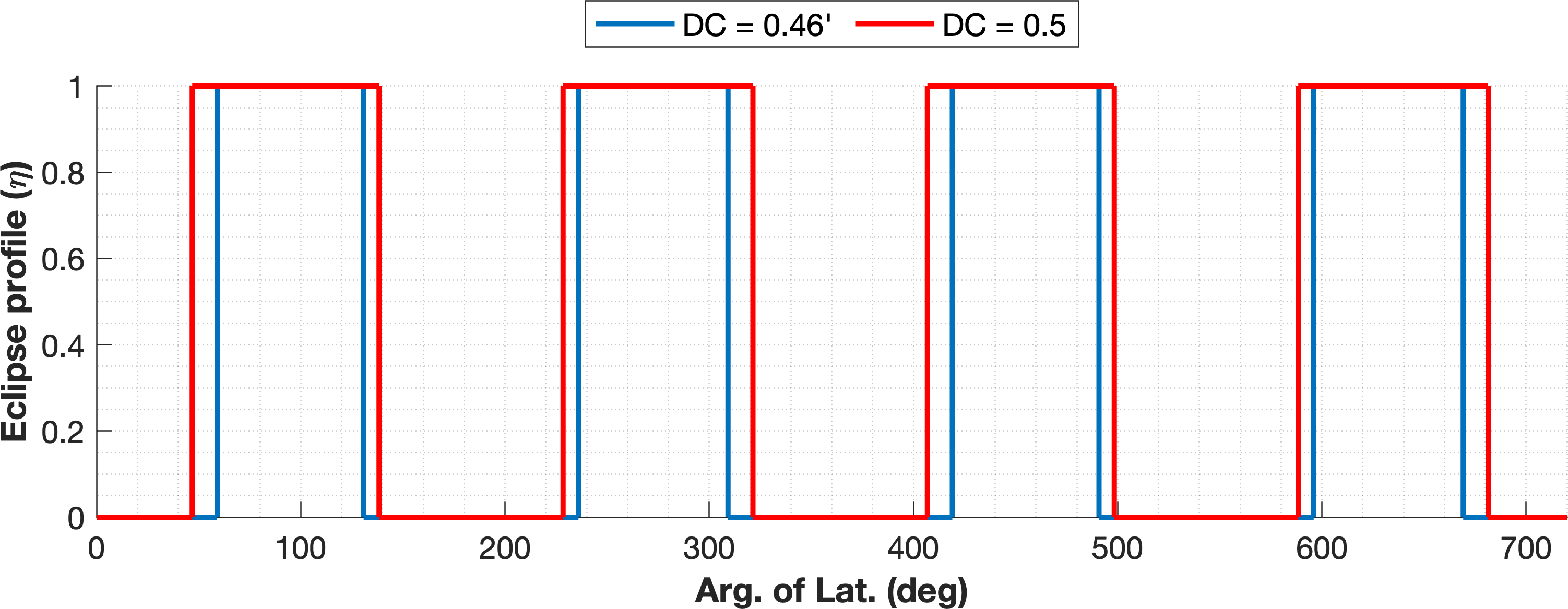}
    \caption{The evolution of $\Delta v'$, semi-major axis, inclination, and RAAN errors with respect to the reference trajectories over time}
    \label{issue2}
  \end{subfigure}
  \caption{Results of the up leg fuel optimal implementation with low thrust errors}
  \label{fig:main2djoij}
\end{figure}

At the start of the MPC guidance, the Extended Edelbaum method and the RAAN matching scheme (collectively called the Preliminary Mission Design Tool (PMDT)) discussed in Algorithms 1 and 2 in [\citenum{ADRMW}] are utilized to obtain a time/fuel optimal trajectory to go from  $a_0, i_0, \Omega_0$ to $a_f, i_f, \Omega_f$. Note that in the RAAN matching process, for certain drift orbit parameters, it may be necessary to wait for many orbits to match RAAN. As such, the search space in terms of TOF has discontinuities, as shown in Figure \ref{issue1}. This may cause convergence problems for gradient-based methods such as fmincon.  Hence, an initial guess for the drift orbit was first obtained through a grid search that finds a minimum $\Delta v/TOF$ that is away from the discontinuity. This was done to direct the fmincon search into an adequate search range.   \par 
The effect of eclipses is incorporated into the PMDT by adjusting the total maximum thrust. However, it does not consider where in the orbit the thrust will be turned off due to eclipses. In reality, if the optimal locations to make orbital changes are in eclipse, making corrections are much harder. To combat this issue, the duty cycle for the reference calculation  ($DC'$) was set to be lower than the real duty cycle ($DC$) of the spacecraft (as shown in Figure \ref{issue2}). While this results in a longer flight time, it creates a positive margin between the reference thrust and the real maximum. This was shown to mitigate the impact of the thrust approximations of the PMDT. 

% When forward propagating the reference PMDT, the effect of eclipses and duty cycle must be considered adequately. However, turning thrust off asymmetrically (i.e., only during the eclipse) will cause eccentricity to build up \cite{Viavattene2022DesignOM}. Hence, thrust is turned off symmetrically across the orbit in the highlighted regions in Figure \ref{eclispeDC}. A detailed description of the eclipse/duty cycle profile allocation can be found in [\citenum{ADRMW}].  

% \begin{figure}[hbt!]
%     \centering
%     \includegraphics[width = 0.5\textwidth]{Figures/eclprop.eps}
%     \caption{Eclipse formulation for propagating with guidance ($C$ denotes the centre of the eclipse, the highlighted regions are where the thrust is turned off.)}
%     \label{eclispeDC}
% \end{figure}

\par 
Finally, it was noted that when the split Edelbaum solution is forward propagated with the thrust it provides, the trajectory undershoots the target. As such, it can be speculated that the $\Delta v$ it estimates may not be sufficient to reach the target in reality. Hence, an additional margin was introduced to the thrust acceleration profile ($f_t^{Edel}$) based on the deviation of the forward propagated PMDT trajectory, as shown in Algorithm \ref{alg:2}.  Note that the time vector for the forward propagation was set such that
$T  = t_0 : P/N: t_f$,
where $P$ and $N$ represent the orbital period at the start and the number of nodes per orbit set for the convex optimization, $t_0$ and $t_f$ represent the start time and the end time of the PMDT reference. 
\par 
\begin{algorithm}[hbt!]
\caption{Forward propagation of the PMDT guess and adjustment of the PMDT thrust profile}\label{alg:2}
\begin{algorithmic}
\Require The PMDT reference (profile of the $t^{Edel}$ and the out of plane angles $\beta^{Edel}$ )
\State Integrate $[\dot{\boldsymbol{X}}, \dot{m}] = fprop(t,x,m)$  to obtain $[\boldsymbol{X},m]$  the osculating orbital elements and mass, using a propagator such as MATLAB's \textit{ode45}.
\State     \hspace{\algorithmicindent}  $[\dot{\boldsymbol{X}}, \dot{m}] = fprop(t,x,m)$:
\State  \hspace{\algorithmicindent}   \hspace{\algorithmicindent}   Calculate the mean keplerian coordinates  $\bar{\boldsymbol{\oe}} = oscCart2meanKep(\boldsymbol{X})$.
\State   \hspace{\algorithmicindent}   \hspace{\algorithmicindent}   Calculate the mean argument of latitude: $\bar{L} =  \bar{\boldsymbol{\oe}}(5) + \bar{\boldsymbol{\oe}}(6) $
\State  \hspace{\algorithmicindent}   \hspace{\algorithmicindent}  Calculate the eclipse/DC profile.
\State \hspace{\algorithmicindent}   \hspace{\algorithmicindent}  \hspace{\algorithmicindent}  Calculate $q_1 = \cos^{-1}(\cos(\bar{L} - L_c)) $  and $q_2 =  \cos^{-1}(\cos(\bar{L}- L_c -\pi )) $ \par \Comment{$L_c$ is the centre of the eclipse, calculated as discussed in [\citenum{ADRMW}]}
 \State \hspace{\algorithmicindent}    \hspace{\algorithmicindent}  \hspace{\algorithmicindent}  \textbf{if} { $q_1 < \frac{\pi}{2}(1-DC')$ \textbf{or}  $q_2 < \frac{\pi}{2}(1-DC')$}  \textbf{then} $\eta=0$ 
\State  \hspace{\algorithmicindent}   \hspace{\algorithmicindent}  \hspace{\algorithmicindent}  \textbf{else} $\eta=1$ \Comment{Note that $\eta=0$ indicates spacecraft being in eclipse.}
\State \hspace{\algorithmicindent}  \hspace{\algorithmicindent}   Compute thrust out of plane angle: $\beta = interp(t^{Edel}, \beta^{Edel}, t_{fprop}(i))$  
 \State  \hspace{\algorithmicindent}   \hspace{\algorithmicindent}  \hspace{\algorithmicindent}  \textbf{if}  { $i_f > i_0$}
 \State  \hspace{\algorithmicindent}   \hspace{\algorithmicindent} \hspace{\algorithmicindent}  \hspace{\algorithmicindent}   \textbf{if}  { $\theta > \pi/2$ or $\theta \leq 3\pi/2$} \textbf{then} $\beta = -|\beta| $  \textbf{else}  { $\beta = |\beta|$}
 \State  \hspace{\algorithmicindent}   \hspace{\algorithmicindent}  \hspace{\algorithmicindent}  \textbf{else}
 \State  \hspace{\algorithmicindent}   \hspace{\algorithmicindent} \hspace{\algorithmicindent}  \hspace{\algorithmicindent}   \textbf{if} { $\theta \geq \pi/2$ or $\theta < 3\pi/2$} \textbf{then}
$\beta = |\beta|$  \textbf{else} 
$\beta = -|\beta|$
\State \hspace{\algorithmicindent}  \hspace{\algorithmicindent}  Compute thrust: $\boldsymbol{a}_T = \eta T_{max}/m[0, \cos{\beta}, \sin{\beta}]$
\State \hspace{\algorithmicindent}  \hspace{\algorithmicindent}  Compute the effect of $J_2$ and drag to the acceleration ($\boldsymbol{a}_{J2}$ and $\boldsymbol{a}_{d}$) as given in [\citenum{8}]. 
\State \hspace{\algorithmicindent}  \hspace{\algorithmicindent}  Calculate the state derivative. $\dot{\boldsymbol{X}}  =   \boldsymbol{a}_T+ \boldsymbol{a}_{J2} + \boldsymbol{a}_{d} $.
\State \hspace{\algorithmicindent}  \hspace{\algorithmicindent}  Determine the derivative of the mass by $\dot{m} = \frac{\eta T_{max}}{I_{sp} g_0}$
\State Calculate the $\Delta v_r = |\bm{\Delta v'}|$ that corresponds to an estimate of the $\Delta v$ required to go from the reached state $\bm{X}(\text{end})$ to the target $a_f, i_f,\Omega_f$ using Equation \ref{dvp}.
\State Compute the adjustment to the Edelbaum $\Delta v$ (${\Delta v^{Edel}}$) as:  $ \Delta v_{adj}(t) = \Delta v^{Edel}(t) + \frac{t}{t_f-t_0}\Delta v_r $
 \State Calculate the adjusted spacecraft mass using the rocket equation:  $m_{adj}(t) = m_0/{\exp(\frac{ \Delta v_{adj}(t) }{I_{sp}g_0})}$
 \State Calculate the adjusted thrust profile:  $f_t^{Edel}(t) = T_{max}/m_{adj}(t)$.
\end{algorithmic}
\end{algorithm}

Note that the exact times when the thrust is turned on and off ($t_{switch}$) are calculated using an event function in ode45 and included in the time vector as
$ T = \textrm{sort}([T, t_{switch}])$. When the thrust profile is discretized in the MPC process without the exact switching times, some thrust will be excluded from the profile, leading to additional tracking difficulties due to insufficient fuel. 

% computing at a lower DR 
\subsection{Convex-based adaptive tracking}
Following the reference generation, convex-optimization-based adaptive tracking is performed. In this section, the trajectory is modeled and the control is optimized for a segment of the total trajectory, and the spacecraft is propagated for this segment of time using the obtained control. The process is then repeated for another segment, and so on. \par
\begin{figure}[hbt!]
    \centering
    \includegraphics[width = \textwidth]{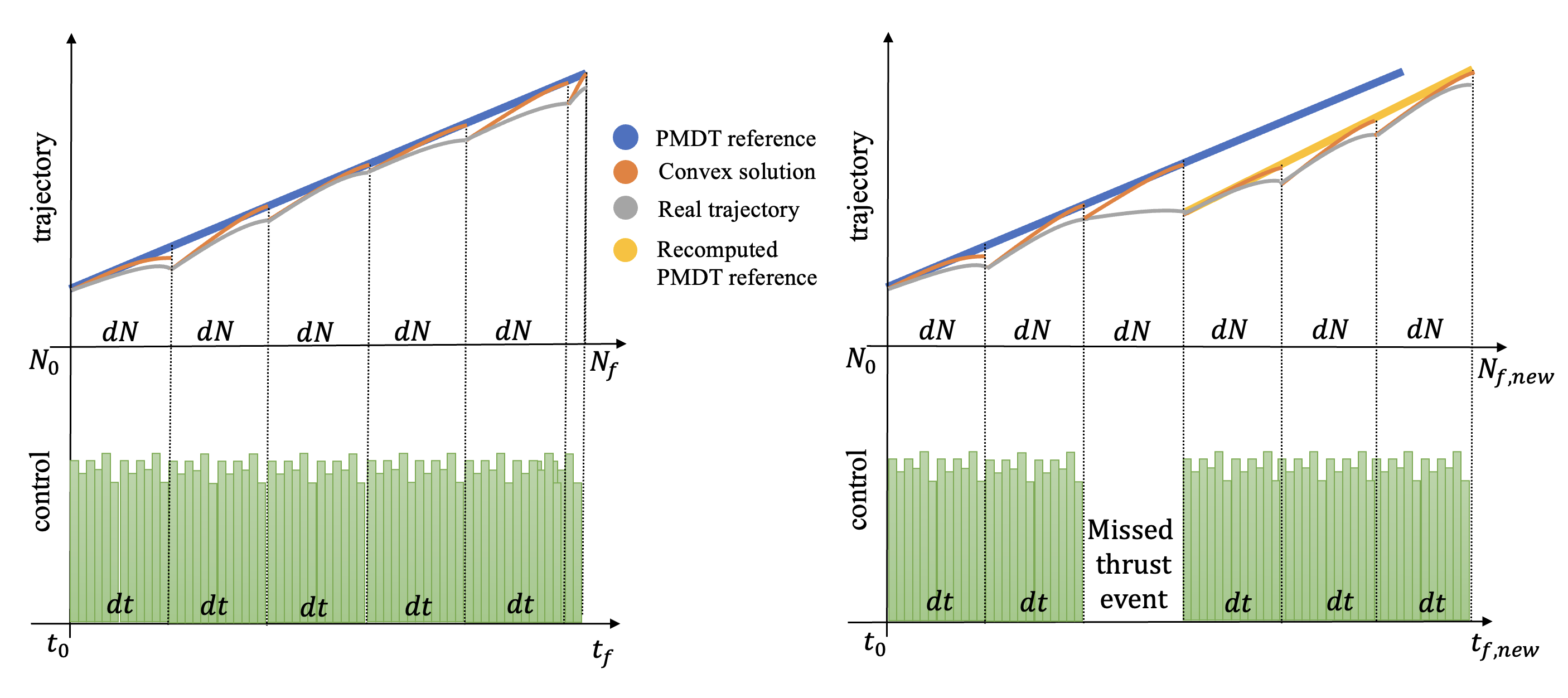}
    \caption{Overview of the convex-based adaptive tracking (left: without reference recomputations, right: with reference recomputation following a misthrust event)}
    \label{MPCprocess}
\end{figure}

\begin{figure}[hbt!]
    \centering

\begin{tikzpicture}[->,>=stealth,shorten >=1pt,auto,node distance=2.5cm,
                    thick,main node/.style={rectangle,draw,font=\sffamily\Large\bfseries}]
      
  \node [process, fill=white!20] at (0,-2) (step1) { \textbf{Setup}: Set the size of the segment ($dN$), $\bm{x}_{0,osc}, m_0$ and set index $N_{run} = 1$};
  
  \node [process, fill=white!20] at (0,-3.5) (step2) {Repeat following while   $N_{run} + dN < \text{len}(T)$: };
   \node[draw, rounded corners] at (0, -6.5) (step3) {
    \begin{minipage}{15cm}
      \begin{myBullets}
      \item  \textbf{Initial guess generation} by propagating $\bm{x}_{0,osc}$ and $m_0$ from time $T(N_{run})$ to $T(N_{run} +dN)$ under $f_t^{Edel}$  (Algorithm \ref{alg:3})
      \item \textbf{Formulate and solve the convex optimization} to obtain the thrust profile ($\bm{a}^{Convex}$) that minimized fuel consumption and the distance to reference at T($N_{run} + dN$).
      \item \textbf {Forward propagation of the spacecraft} under $\bm{a}_{Convex}$ to simulate motion (Algorithm \ref{alg:6}).  Obtain the reached osculating state, $\bm{x_{f,osc}}$, mass at the end ($m_f$) and the approximate $\Delta v'$ required to go from $\bm{x_{f,osc}}$ to the target state for the segment (calculated using Equation \ref{dvp}.). 
      \end{myBullets}
    \end{minipage}
  }
  ;
   \node [process, fill=white!20] at (-5,-10) (dvgood) { $\bm{x}_{0,osc} = x_{f,osc}$, $m_0 = m_f$,  $N_{run} = N_{run} + dN$ };
      \node [process, fill=white!20] at (4,-10) (dvbad) { \textbf{Reference regeneration} using PMDT and Algorithm \ref{alg:2} };
  %  \node [process, fill=white!20] at (0,-2) (dvgood) { $\bm{x}_{0,osc} = x_{f,osc}$, $m_0 = m_f$,  $N_{run} = N_{run} + dN$ };
    \node[draw=none] at (-2.5, -9) (no) {$\Delta v' < \epsilon$ };
\node[draw=none] at (4, -9) (yes) {$\Delta v' > \epsilon$};

         \node[draw=none] at (-9, -10.15) (fakenode) {};
         
\path [connector] (step1) -- (step2);
\path [connector] (step2) -- (step3);
\path [connector] (step3) -- (dvgood);
\path [connector] (step3) -- (dvbad);

\draw[-] (dvgood.west) -- (fakenode.north);

%\path [connector] (dvgood.west) -- (fakenode.north);

\path [connector] (fakenode) |- (step2);
\end{tikzpicture}
\caption{Convex-based adaptive tracking algorithm}
\label{2}
\end{figure}
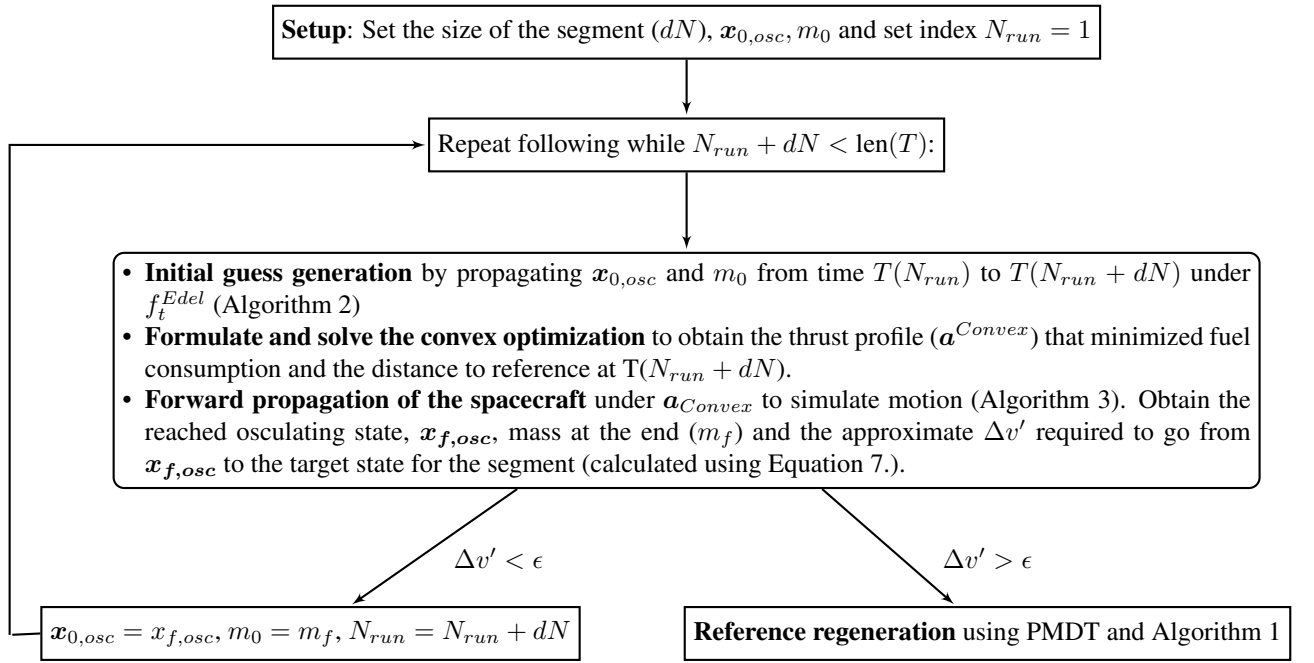

 
Each segment includes n orbits, and each orbit is discretized into N nodes. As such, each segment is of grid size $dN = nN$.   Convex optimization is used to optimize the profile of control accelerations ($\bm{a}^{convex}$) within the segment. Then the spacecraft is forward propagated to simulate its motion during the optimized segment with the convex-optimized control accelerations. In situations where the spacecraft has deviated from the reference trajectory significantly, the reference PMDT is recomputed such that the spacecraft can still reach the target at the cost of increased TOF and/or $\Delta v$, as shown in Figure \ref{MPCprocess}. Note that the $\Delta v'$ parameter estimates the tracking error and determines when a recomputation of the reference trajectory will be needed. 
Figure \ref{2} shows the overview of the guidance process.

\definecolor{BulletsColor}{rgb}{0, 0, 0}
\newlist{myBullets}{enumerate}{1}

\setlist[myBullets]{
  label=\textcolor{BulletsColor}{\textbullet},
  leftmargin=*,
  topsep=0ex,
  partopsep=0ex,
  parsep=0ex,
  itemsep=0ex,
  before={\color{BulletsColor}}
}

\textit{Setup} \par 
Before the adaptive tracking is run, Let $ dN = nN$. 
The starting osculating state is set as $\boldsymbol{x}_{0,osc}$ and the starting mass ($m_0$) is set as the initial spacecraft mass.

\textit{Initial guess generation} \par 
It was seen that the eclipse profile from the forward propagation of PMDT from earlier is not sufficiently accurate for the convex optimization, especially when the real trajectory starts to deviate from the reference.  Hence, for each segment, an initial guess is generated by forward propagating the thrust accelerations provided by the PMDT reference, starting from $\bm{x}_{0,osc}$. The initial guess generation for one segment is discussed in Algorithm \ref{alg:3}, where step-wise integration with a zero-order hold on the thrust acceleration is utilized. This generates an initial guess that matches the PMDT reference and adheres to the dynamics even when  $\bm{x}_{0,osc}$ differs from the reference, hence obtaining a more accurate profile of eclipses ($\bm{\eta}_{guess}$).
\begin{algorithm}[hbt!]
\caption{Initial guess generation for a single run of adaptive tracking}\label{alg:3}
\setstretch{0.75}
\begin{algorithmic}
\Require starting coordinates ($\bm{x}_s = [\bm{r}_s, \bm{v}_s]^T = \boldsymbol{x}_{0,osc}$) and starting mass ($m_0$)
\State Define the time vector for this segment. $t = T(N_{run} : N_{run} + dN) $.
\State Let $m(1) = m_0$. 
\State Calculate the mean Keplerian coordinates  $\bar{\boldsymbol{\oe}} = oscCart2meanKep(\boldsymbol{x}_{s})$.
\For{$i = 1: dN$}
\State  Calculate the mean argument of latitude: $\bar{L} =  \bar{\boldsymbol{\oe}}(5) + \bar{\boldsymbol{\oe}}(6) $
\State   Calculate the eclipse/duty cycle profile.
\State \hspace{\algorithmicindent} Calculate $q_1 = \cos^{-1}(\cos(\bar{L} - L_c)) $  and $q_2 =  \cos^{-1}(\cos( \bar{L}- L_c -\pi )) $ \par \Comment{$L_c$ is the centre of the eclipse, calculated as discussed in [\citenum{ADRMW}]}
 \State  \hspace{\algorithmicindent}  \textbf{if} { $q_1 < \frac{\pi}{2}(1-DC')$ \textbf{or}  $q_2 < \frac{\pi}{2}(1-DC')$}  \textbf{then} $\eta=0$ 
 \State   \hspace{\algorithmicindent}  \textbf{else} $\eta=1$ \Comment{Note that $\eta=0$ indicates spacecraft being in eclipse.}
\State Using the thrust acceleration profile from the Edelbaum guess, calculate the acceleration $ a(i)$.
\begin{equation}
f_s (i)  = \textrm{interp}(t^{Edel} , f^{Edel}_t , t(i)) \ \text{and} \ f_s (i+1)  = interp(t^{Edel} , f^{Edel}_t , t(i+1))
\end{equation}
\begin{align} a(i) = \frac{1}{DC'} \eta(i) \frac{f_s(i+1) + f_s(i)}{2} \end{align}
\State Using the $\Delta v(i) = a(i) (t(i+1) -t(i))$ and the rocket equation, determine the mass after the burn: $m(i+1) = \frac{m(i)}{\exp(\Delta v/I_{sp}/g_0)}$
\State Determine the average out-of-plane angle for the burn. 
\State \hspace{\algorithmicindent}  Calculate the out of plane angles at $t(i)$ and $t(i+1)$ (called $\beta(i)$ and $\beta(i+1)$): 
 {
\begin{align} 
\beta(i) &= \tan^{-1} \left( \frac{V_0 \sin{\beta_0}}{V_0 \cos{\beta_0} - \frac{T_{max}}{m(i)}t(i)}\right) 
\end{align}}
\State \hspace{\algorithmicindent} Calculate the average out-of-plane angle as: 
$\beta_{avg} = \frac{1}{2}(\beta(i)  + \beta(i+1))$
\State Determine the direction of the out-of-plane angle.

\State \hspace{\algorithmicindent} \textbf{if} { $i_f > i_0$} \textbf{then}
\State \hspace{\algorithmicindent} \hspace{\algorithmicindent}  \textbf{if} { $\bar{L} > \pi/2$ or $\bar{L}  \leq 3\pi/2$} \textbf{then} $\beta_{avg} = -|\beta_{avg}|$ \textbf{else}{ $\beta_{avg} = |\beta_{avg}|$}
\State \hspace{\algorithmicindent} \textbf{else}
\State \hspace{\algorithmicindent}  \hspace{\algorithmicindent}  \textbf{if} { $\bar{L}  \geq \pi/2$ or $\bar{L}  < 3\pi/2$ } \textbf{then}
$\beta_{avg} = |\beta_{avg}|$  \textbf{else} { $\beta_{avg} = -|\beta_{avg}|$}

\State Calculate the $a$ components in RTN:  { $\boldsymbol{a}_{RTN} = a(i) [0 , \cos{\beta_{avg}}, \sin{\beta_{avg}}] $
\State Let $\boldsymbol{x}_{guess}(i,:) = [\boldsymbol{x}_s, \boldsymbol{a}_{RTN} ]$}
\State  {Propagate from $t(i)$ to $t(i+1)$ under $J_2$, drag and gravitational acceleration, converting the RTN accelerations to ECI.
$\bm{a}_{ECI} = [\frac{\bm{{r}_s}}{|\bm{{r}_s}|} , \frac{\bm{h}_s \times \bm{r}_s}{|\bm{h}_s \times \bm{r}_s|},  \frac{\bm{{h}_s}}{|\bm{{h}_s}|} ]\bm{a}_{RTN}$
}
\State  { Obtain the state at $t(i+1)$ and define it as the new $\boldsymbol{x}_s$.} 
\State Obtain the new $\bar{\boldsymbol{\oe}}$ by converting $\boldsymbol{x}_s$ to a mean keplerian state.
\EndFor
\State Let $\boldsymbol{x}_{guess}(\text{end},:) = [\boldsymbol{x}_s, 0, 0, 0]$ \Comment{No thrust is allocated for the last node.}
\State \textbf{Output:} $\boldsymbol{x}_{guess}$.
\end{algorithmic}
\end{algorithm}

\textit{Formulating and solving the convex optimization problem} \par 
Once a segment's initial guess is obtained, convex optimization is used to obtain the optimal acceleration control that minimizes the total fuel consumption while remaining close to the PMDT reference. As such, the objective function is defined as \par 
\begin{equation}
    \textrm{minimize } J = \int^{t_s +dt}_{t_s} 	\lVert \bm{a}\rVert dt  + \lVert \bm{\Delta v'} \rVert 
\end{equation}
subject to
\begin{equation}
    \begin{aligned}
& \dot{\bm{x}}=\boldsymbol{f}(\boldsymbol{x}, \boldsymbol{u})= \bm{h(x)} + \bm{a} \\
& \bm{x} \left(t_s\right)= \bm{x}_{0,osc} \\
& \|\boldsymbol{a}\| \leq \frac{T_{max}}{m} \\
\end{aligned}
\end{equation}
where $t_s$ is the starting time step of the segment and $dt$ is the time length of the segment.  $\bm{h(x)}$ is the natural dynamics. which is defined as $\bm{h(x)} = \bm{a}_d + \bm{a}_{J2}$ as given in Algorithm 1. $\bm{a} = [a_r, a_\theta, a_N]^T$ and is the thrust acceleration in RTN coordinates. \par 
 Note that the convexification is done around the initial guess obtained in the previous step. For each segment, the target coordinates were calculated as:  
$$\oe_t = \textrm{interp}(t^{Edel}, \oe^{Edel}, t_s + dt)$$
\par 
The convex tracking is only required to match the reference $a, i$ and $\Omega$ precisely once the spacecraft reaches its target. Precise matching to the PMDT trajectory before reaching the target would result in unnecessary fuel consumption. Hence, the margin of $\Delta v$ between the state reached at the end of the convex optimization, and $\oe_t$ is calculated and minimized in the convex optimization. This quantity is labeled as $\Delta v'$ and is calculated by integrating the modified equinoctial Gauss Variational Equations of maximum rates of change \cite{GVE} in time, %as follows, assuming the states are time-independent.
%     \begin{align}
%         \dot{a}& =2 f_s a_t \sqrt{\frac{a_t}{\mu}} \sqrt{\frac{1+e_t}{1-e_t}} &\rightarrow  \Delta a = 2 \Delta v_{a_t} a_t \sqrt{\frac{a_t}{\mu}} \sqrt{\frac{1+e_t}{1-e_t}}  \\ 
% \dot{h}& =\frac{1}{2} f_s \sqrt{\frac{p_t}{\mu}} \frac{s_t^2}{\sqrt{1-g_t^2}+f_t} &\rightarrow \Delta h_t  =\frac{1}{2} \Delta v_h \sqrt{\frac{p_t}{\mu}} \frac{s_t^2}{\sqrt{1-g_t^2}+f_t}  \\
% \dot{k} & =\frac{1}{2} f_s \sqrt{\frac{p_t}{\mu}} \frac{s_t^2}{\sqrt{1-f_t^2}+g_t} &\rightarrow  \Delta k = \frac{1}{2} \Delta v_k \sqrt{\frac{p_t}{\mu}} \frac{s_t^2}{\sqrt{1-f_t^2}+g_t} 
%     \end{align}
% Then, $\bm{\Delta v'}$ is calculated as: 
providing:
    \begin{align}
     \bm{\Delta v'} &= [{\Delta v_a}, {\Delta v_h},{\Delta v_k}]^T\\  &= \left[ \frac{\Delta a}{2 a }\sqrt{\frac{\mu}{a}}\sqrt{\frac{1-e}{1+e}},{2\Delta h}\sqrt{\frac{\mu}{p}} \frac{\sqrt{1-g^2}+f}{s^2} ,{2\Delta k}\sqrt{\frac{\mu}{p}} \frac{\sqrt{1-f^2}+f}{g^2} \right]^T
     \label{dvp}
 \end{align}
      
% Note that the assumption of time-independent states is valid only over small periods of time.
$ \bm{\Delta v'}$  provides a rough measure of the required $\Delta v$  to reach the reference end state for each segment. Furthermore, it provides a single parameter of accuracy, which is used to determine when a recomputation of the reference is needed in the MPC guidance process. \par
Now that the problem has been defined, it must be convexified to utilize convex optimization in the solution process. The problem dynamics can be convexified as 

\begin{equation}
\bm{x} ({i+1}) = \bm{A} (i)(\bm{\hat{\bm{x}}},\hat{\bm{a}})\bm{x} (i) + \bm{B} (i)(\bm{\hat{\bm{x}}},\hat{\bm{a}}) {\bm{a} (i)}  
\end{equation}
where $\bm{A} (i)(\bm{\hat{\bm{x}}},\hat{\bm{a}})$ and $\bm{B} (i)(\bm{\hat{\bm{x}}},\hat{\bm{a}})$ are the state transition matrices evaluated at the initial guess trajectory, $[\bm{\hat{\bm{x}}},\hat{\bm{a}}]$.
The thrust acceleration constraint can be transformed into a convex cone constraint as 
\begin{equation}
\begin{gathered}
 \sqrt{a_R(i) ^2 + a_T(i) ^2 + a_N(i) ^2}  \leq \widetilde{a}(i) \\
0 \leq \widetilde{a}(i)  \leq T_{max}/m(i) 
\end{gathered}
\end{equation}
where the mass evolution $m$ is obtained from the initial guess. \par 
The $\Delta v'$ constraint shall also be written as a convex cone

\begin{equation}
     \sqrt{\Delta v_a^2 + \Delta v_h^2 + \Delta v_k^2}  \leq \widetilde{\Delta v'}. \\
\end{equation}
Then the convex problem becomes 
\begin{equation}
   \text{Optimize} \   K = \sum_{i=1}^{dN} \widetilde{a}(i) {dt}(i) +  \widetilde{\Delta v'} \ \text{where }  dt(i) = t(i+1) - t(i)
\end{equation}
subject to
\begin{equation}
    \begin{aligned}
&\bm{x} ({i+1}) = \bm{A} (i)(\bm{\hat{\bm{x}}},\hat{\bm{a}})\bm{x} (i) + \bm{B} (i)(\bm{\hat{\bm{x}}},\hat{\bm{a}}) {\bm{a} (i)}  {\bm{a}_i}   \ \text{for } i = 1,...,dN \\
& \bm{x}(1)= \bm{x}_{0,osc} \\
& \sqrt{a_R(i)^2 + a_T(i)^2 + a_N(i)^2}  \leq \widetilde{a(i)} \ \text{for } i = 1,...,dN  \\
& 0 \leq \widetilde{a}(i) \leq T_{max}/m(i) \\
& \sqrt{\Delta v_a^2 + \Delta v_h^2 + \Delta v_k^2}  \leq \widetilde{\Delta v'} \\
\end{aligned}
\end{equation}

Following the convexification, a convex solver such as MOSEK \cite{mosek} or Gurobi \cite{gurobi} can be utilized to solve it. In this work, the state transition matrices involved were obtained using DACE \cite{dace}, and MOSEK was used  perform convex optimization.

\textit{Avoiding successive convex iterations} \par 
As mentioned, the nonconvex nature of the problem dynamics and constraints makes convex optimization unreliable. As such consecutive iterations are required to converge to a near-optimal solution. However, using successive iterations is not ideal for real-time guidance, as it is time-consuming and convergence is not guaranteed.  \par 
Successive iterations are avoided in this work by using the highly accurate  PMDT reference as an initial guess and utilizing a near-linear set of coordinates to represent the problem dynamics. Specifically, the Generalized Equinoctial Orbital Elements (GEqOE) introduced in [\citenum{claudio}] are used as they provide an almost-linear representation of the trajectory even in the presence of $J_2$ perturbations.  \par 
By definition, the GEqOE coordinates are: 
\begin{equation}
    \begin{array}{lll}
\nu, & p_1-g \sin \Psi, & p_2-g \cos \Psi, \\
\mathcal{L}, & q_1-\tan \frac{i}{2} \sin \Omega, & q_2-\tan \frac{i}{2} \cos \Omega .
\end{array}
\end{equation}
where $g, \nu,\Psi, \mathcal{L}$ are the eccentricity vector, generalized mean motion, generalized longitude of periapsis, and generalized mean longitude. The effect of $J_2$ and other perturbations are included in the dynamics through formulating them as contributions to the total potential energy \cite{claudio}. As such, the generalized mean motion state variable remains constant even in the presence of $J_2$ \cite{claudio2}. \par
\begin{figure}[hbt!]
            	\centering
            	\includegraphics[width=0.6\textwidth]{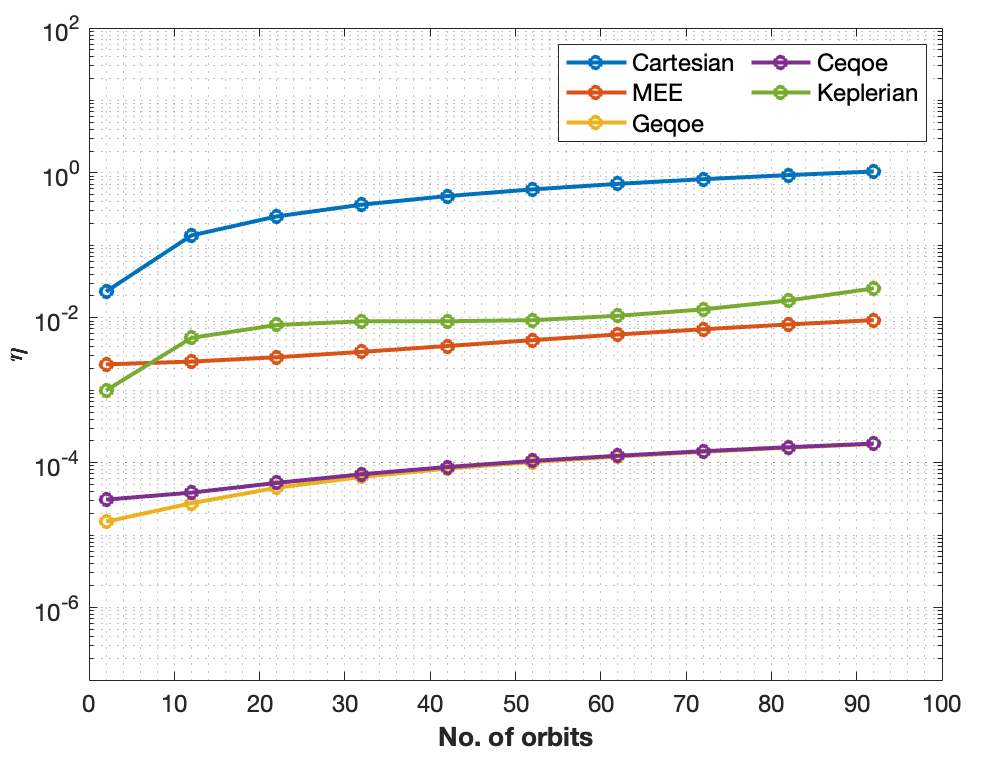}
            	\caption{Comparison of the nonlinearity index (as defined in [\citenum{Junkins2003HowNI}])  of GEqOE with other coordinates. }
             \label{Geqfig}
\end{figure}

 
To further confirm the near-linear nature of GEqOE and how it compares to other coordinate systems, its nonlinearity index \cite{Junkins2003HowNI} was compared to other coordinate systems for a J2-only perturbed two-body system. A propagation was done for a spacecraft at an altitude of 350 km with an inclination of 99.22 deg in Cartesian, Modified Equinoctial, Generalised Equinoctial, Classical Equinoctial and Keplerian coordinates.  Initial normalized position and velocity uncertainties of 1 km and 1 m/s were then introduced to analyze how the uncertainty is propagated forward in each coordinate system.  The propagated uncertainty of each of the coordinates was estimated using the nonlinearity index ($v$) defined by J. Junkins in [\citenum{Junkins2003HowNI}], where $v$ was defined as  

\begin{equation}
    v \triangleq \sup _{i=1, \cdots, N} \frac{\left\|A\left(\mathbf{x}_{\mathrm{i}}\right)-A(\overline{\mathbf{x}})\right\|}{\|A(\overline{\mathbf{x}})\|}
\end{equation}
where $A$, $\mathbf{x}$ and $\overline{\mathbf{x}}$ are the state transition matrix, the perturbed state, and the unperturbed state after propagating for 15 orbits. The results are shown in Figure \ref{Geqfig}, which illustrates that the nonlinearity index of GEqOE is much lower than other coordinate systems, especially for a lower number of orbits propagated.  As such, using it in the context of convex optimization enhances the validity of the obtained solution by reducing the nonlinearity involved.

% The Generalised Equinoctial Orbital Elements (GEqOE) \cite{claudio} are used in the convex formulation to retain the dynamics in a near-linear format (as shown in Figure \ref{Geqfig}), enhancing the validity of the convex solution compared to other coordinate systems.  

\textit{Forward propagation of the spacecraft} \par 
Once the optimal acceleration profile is obtained using convex optimization, the spacecraft is propagated under realistic dynamics and thrust uncertainties. In this section, thrust bias and magnitude errors are considered, as well as events of misthrust where the thrust is turned off for a full tracking segment. Algorithm \ref{alg:6} illustrates the process of forward propagation for a single segment. \par 
Note that the spacecraft is separately propagated under a low-fidelity model and a high-fidelity model. The low-fidelity model has the same dynamics as the convex optimization and considers only the forces due to drag, $J_2$ perturbations, and thrust. The high-fidelity model used is the one developed in [\citenum{MORSELLI2014490}], which takes the full geopotential (zonal, tesseral and sectorial harmonics), solar radiation pressure and the third body perturbations of the sun and the moon into account. As the high-fidelity model calculates the geopotential in the ECEF frame,  the effects of precession and nutation are also considered. Hence, a conversion from ECI to TOD (true of date) coordinates must be implemented following the high-fidelity propagation.\par
 
\begin{algorithm}[hbt!]
\caption{Forward propagation of the spacecraft under convex optimized acceleration}\label{alg:6}
\begin{algorithmic}
\State 
\Require $\bm{x}_s =  [\bm{r}_s , \bm{v}_s]^T = \boldsymbol{x}_{0,osc} $, $m_0$ and the convex optimized control profile $\bm{a}^{Convex}$.
\State Let $m(1) = m_0$. 
\State Allow misthrust events as: 
\If{$x \sim \mathcal{U}(0,1) < p_{misthrust}$} 
$\bm{a}^{Convex} = \bm{0}$  \Comment{$ p_{misthrust}$ is the user-set misthrust probability}
\EndIf
\For{$i = 1: dN$}
\State Recompute the real eclipse/Duty cycle ($\eta$) using $DC$ as shown in Algorithm \ref{alg:3}.
\State Introduce normalized thrust magnitude and direction errors.
\State \hspace{\algorithmicindent} Compute the thrust magnitude with error:
\State \hspace{\algorithmicindent} \hspace{\algorithmicindent}  $  |\bm{a}^{Convex} (i)|_{E}   = |\bm{a}^{Convex} (i)|   + |\bm{a}^{Convex} (i) | \delta_T \ \text{where}  \ \delta_T  \sim\mathcal{N}(0,\bar{\delta_T})$  \State\hspace{\algorithmicindent} \hspace{\algorithmicindent}  \Comment{$\bar{\delta_T}$:fractional standard deviation of thrust magnitude error}
\State \hspace{\algorithmicindent} Compute the out-of-plane angle:  $\beta = \sin^{-1}( \frac{\bm{a}^{Convex}(i)_N}{|\bm{a}^{Convex} (i)|})$
\State \hspace{\algorithmicindent} Compute the in-plane angle:  $\alpha = \tan^{-1}( \frac{\bm{a}^{Convex}(i)_T}{\bm{a}^{Convex}(i)_N})$
\State \hspace{\algorithmicindent} Calculate the out-of-plane angle with error: 
\State \hspace{\algorithmicindent}  \hspace{\algorithmicindent}  $\beta_E    = \beta  + \delta_\beta   \ \text{where}  \  \delta_\beta \sim \mathcal{N}(0,\bar{\delta_\beta})$
\State \hspace{\algorithmicindent} \hspace{\algorithmicindent}  \Comment{$\bar{\delta_\beta}$: standard deviation of the out of plane thrust angle error}
\State \hspace{\algorithmicindent} Compute the thrust acceleration with errors: 
\State  \hspace{\algorithmicindent}  \hspace{\algorithmicindent}  $\bm{a}^{Convex} (i)_E = |\bm{a}^{Convex} (i)|_{E}  [\cos{\beta_E}\sin{\alpha}, \sin{\alpha}\cos{\beta_E}, \sin{\beta_E}]^T$
% \State  {Calculate $\Delta v = \eta a (i) (t(i+1) - t(i))$}
% \State Update mass. $ m(i) = \frac{m(i-1)}{\exp(\Delta v/I_{sp}/g_0)}$.
\If {Low Fidelity propagation}
\State  {Propagate from $t(i)$ to $t(i+1)$ under $J_2$, drag, gravitational acceleration and $\eta \bm{a}^{Convex} $ (converted to ECI) to get the state at $t(i+1)$.}  
\EndIf
\If {High Fidelity propagation}
\State Convert $\boldsymbol{x}_{fprop}(i,:)$ from TOD to ECI. 
\State  {Propagate from $t(i)$ to $t(i+1)$ under $\bm{a} $ using the high fidelity dynamics model in [ \citenum{MORSELLI2014490}].}  
\State Convert the output state from ECI to TOD.
\EndIf
\State  Define the new cartesian state at t(i+1) as $\bm{x}_s$.
\EndFor
\State Calculate $\Delta v'$ to go from $\bm{x}_s$ to the target state for this segment.
\State Output $\boldsymbol{x}_{f,osc} = \bm{x}_s, \Delta v'$ and $m_f = m(end)$.

\end{algorithmic}
\end{algorithm}

% The $\Delta v'_{fprop}$ is used to determine whether a reference recomputation is required after the spacecraft has been forward propagated. If $\Delta v'_{fprop} > \epsilon $, a recomputation shall be conducted. 

\subsection{Reference Recomputation}
As the error builds up due to the presence of thrust uncertainties and nonlinear dynamics encountered, the $\Delta v'_{fprop}$ value increases over time. As the reference duty cycle ($DC'$) is lower than the spacecraft's real duty cycle ($DC$), the additional thrust from this margin resolves some trajectory error buildup. However, for long-duration transfers that entail large changes in orbital parameters, a recomputation of the reference mid-trajectory may be essential to keep the spacecraft on a course that meets the target. \par 
As such, a tolerance ($\epsilon$) is set on $\Delta v'$, such that once exceeded, the PMDT reference and its forward propagation are recomputed to ensure that the spacecraft remains in a trajectory that achieves the target semi-major axis, inclination and RAAN. 

\section{Results}

In this section, results are provided for an upward and downward leg of the multi-debris removal tour. The upward leg takes the servicer spacecraft from a 350 km altitude orbit to the orbit of the H-2A (F15) rocket body, and the downward leg brings the debris and the servicer back down to an orbit of 350 km altitude. Thirty days are allocated for proximity and docking operations at the H-2A (F15) orbit between the up and down legs. The simulations were conducted with no thrust errors with two sets of thrust errors, as shown in Table \ref{tError}.

\begin{table}[hbt!]
\centering 
\begin{tabular}{ccccc}
Simulation         & $DC'$  & $p_{misthrust}$ (\%) &  $\bar{\delta_T}$ (\%) & $\bar{\delta_\beta}$(deg) \\ \hline 
No thrust errors    & 0.40 & 0                 & 0                                                    & 0                                                                        \\
Low thrust errors & 0.40 & 1                 & 5                                                    & 5                                                                        \\
High thrust errors & 0.40 & 7                & 7                                                   & 7   \\ \hline                                                                    
\end{tabular}
\caption{Thrust error settings utilized in the case study}
\label{tError}
\end{table}

As the thrust errors are stochastic, three simulations were conducted with the same set of errors for both high and low-fidelity propagations. The results given in Tables \ref{tR} and \ref{td} for high and low thrust errors are the averaged outcomes of these simulations. 
% Monte Carlo simulations shall be conducted on the thrust errors in the future to obtain the relevant statistical information about the propagations with errors. 
In the convex tracking, the number of nodes per orbit (N) was set to 36. The number of orbits per segment (n) was set to 5. Note that $\epsilon$ was set to be 2 m/s to determine when recomputations will be needed.

\subsection{Fuel optimal up leg from the initial orbit to the target debris }

In this section, the up leg of the mission was optimized for fuel consumption using the PMDT, and the obtained trajectory was used as a reference for the convex-MPC guidance. The upwards leg starts at 00:00 UTC on 25 March 2022, and the orbital parameters of the  H-2A (F15) rocket body and the starting orbit are given in Table \ref{T1}. Note that while only semi-major axes, inclination, and RAAN are provided as targets, the eccentricity is maintained to be quasi-circular by alternating the direction of the out-of-plane thrust as shown in Algorithm \ref{alg:3}. The servicer spacecraft is assumed to have a wet mass of $m_0 = $ 800 kg, with a drag coefficient of 2.2 and a frontal area of $\SI{0.01}{\meter\squared}$. The low thrust engine of the spacecraft is assumed to have 60 mN maximum thrust, 1300 s $I_{sp}$, and a duty cycle ($DC$) of 50\%.

\begin{table}[hbt!] \centering\begin{tabular}{ccccccc}Object & $a$ (km) & e & $i$ (deg) & $\Omega$ (deg) & $\omega$ (deg) & $\theta$ (deg) \\ \hline Initial & 6728.1363 & 0.004 & 98.3 & 15.3 & 0 & 0 \\Target(H-2A F-15) & 6975.0874 & 0.0040111 & 98.1521 & 19.9669 & - & - \\\hline \end{tabular}\caption{{Orbital elements of initial orbit and the target debris at launch time}}
\label{T1} \end{table}

\begin{figure}[hbt!]
  \centering
  \begin{subfigure}[b]{0.49\textwidth}
  \centering
    \includegraphics[width = \textwidth]{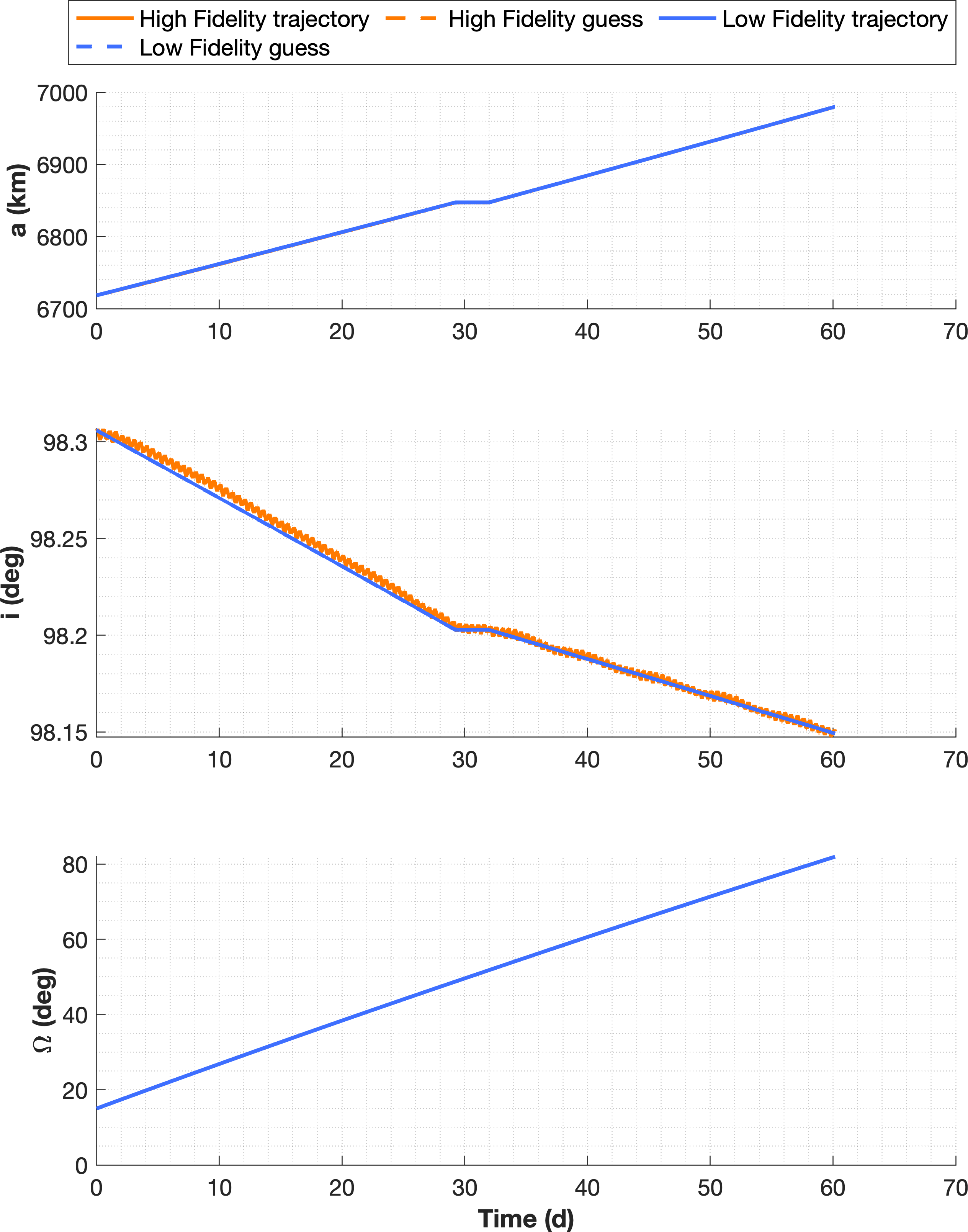}
    \caption{The fuel optimal solution with the high and low fidelity spacecraft propagations}
    \label{foptb1}
  \end{subfigure}
  \hfill
  \begin{subfigure}[b]{0.49\textwidth}
    \centering
    \includegraphics[width = \textwidth]{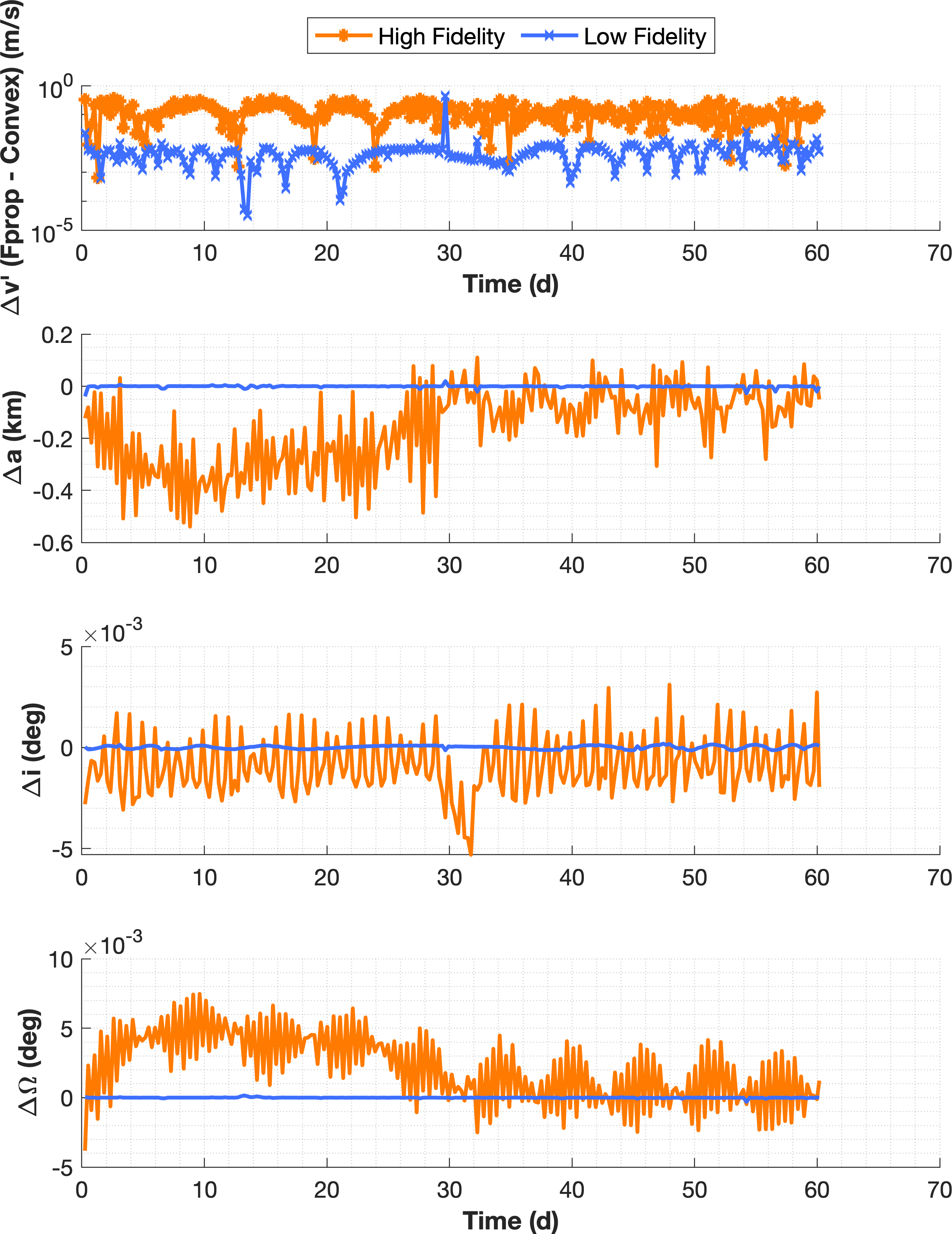}
    \caption{The evolution of $\Delta v'$, semi-major axis, inclination, and RAAN errors with respect to the reference trajectories over time}
    \label{foptb2}
  \end{subfigure}
  \caption{Results of the up leg fuel optimal implementation without thrust errors}
  \label{fig:main1}
\end{figure}

\begin{figure}[hbt!]
  \centering
  \begin{subfigure}[b]{0.49\textwidth}
  \centering
    \includegraphics[width = \textwidth]{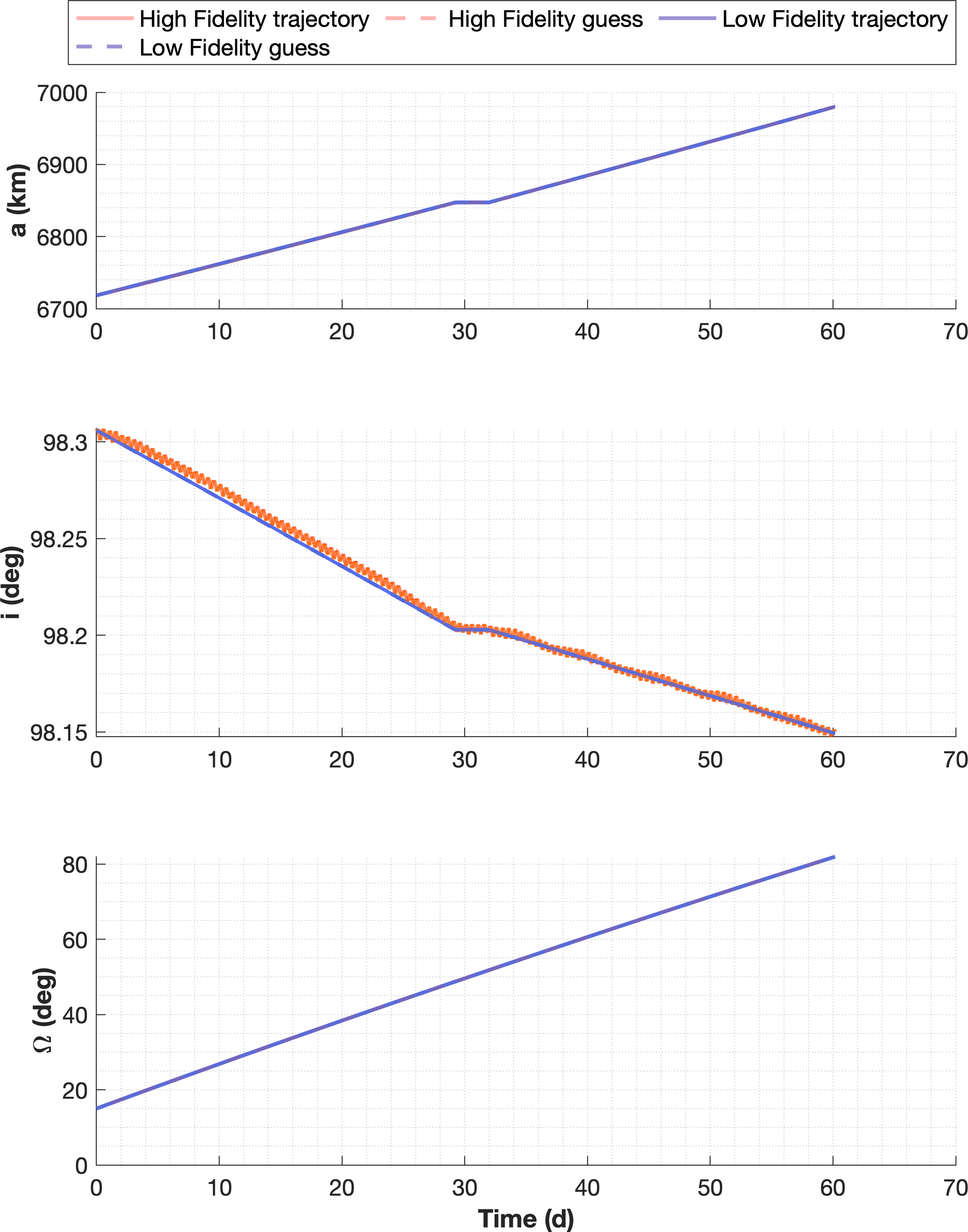}
    \caption{The fuel optimal solution with the high and low fidelity spacecraft propagations}
    \label{foptc11}
  \end{subfigure}
  \hfill
  \begin{subfigure}[b]{0.49\textwidth}
    \centering
    \includegraphics[width = \textwidth]{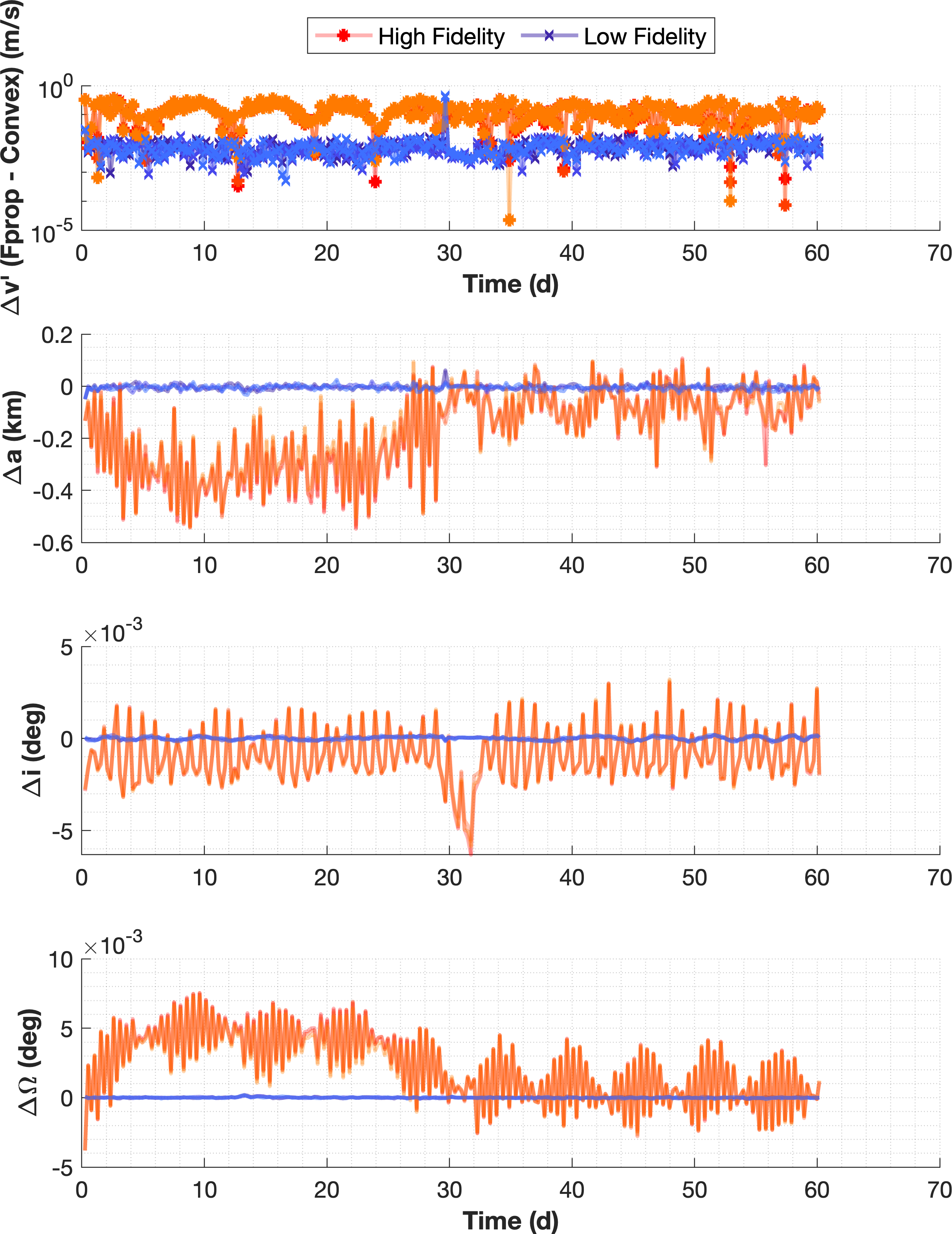}
    \caption{The evolution of $\Delta v'$, semi-major axis, inclination, and RAAN errors with respect to the reference trajectories over time}
    \label{foptc21}
  \end{subfigure}
  \caption{Results of the up leg fuel optimal implementation with low thrust errors}
  \label{fig:main2}
\end{figure}

\begin{figure}[hbt!]
  \centering
  \begin{subfigure}[b]{0.49\textwidth}
  \centering
    \includegraphics[width = \textwidth]{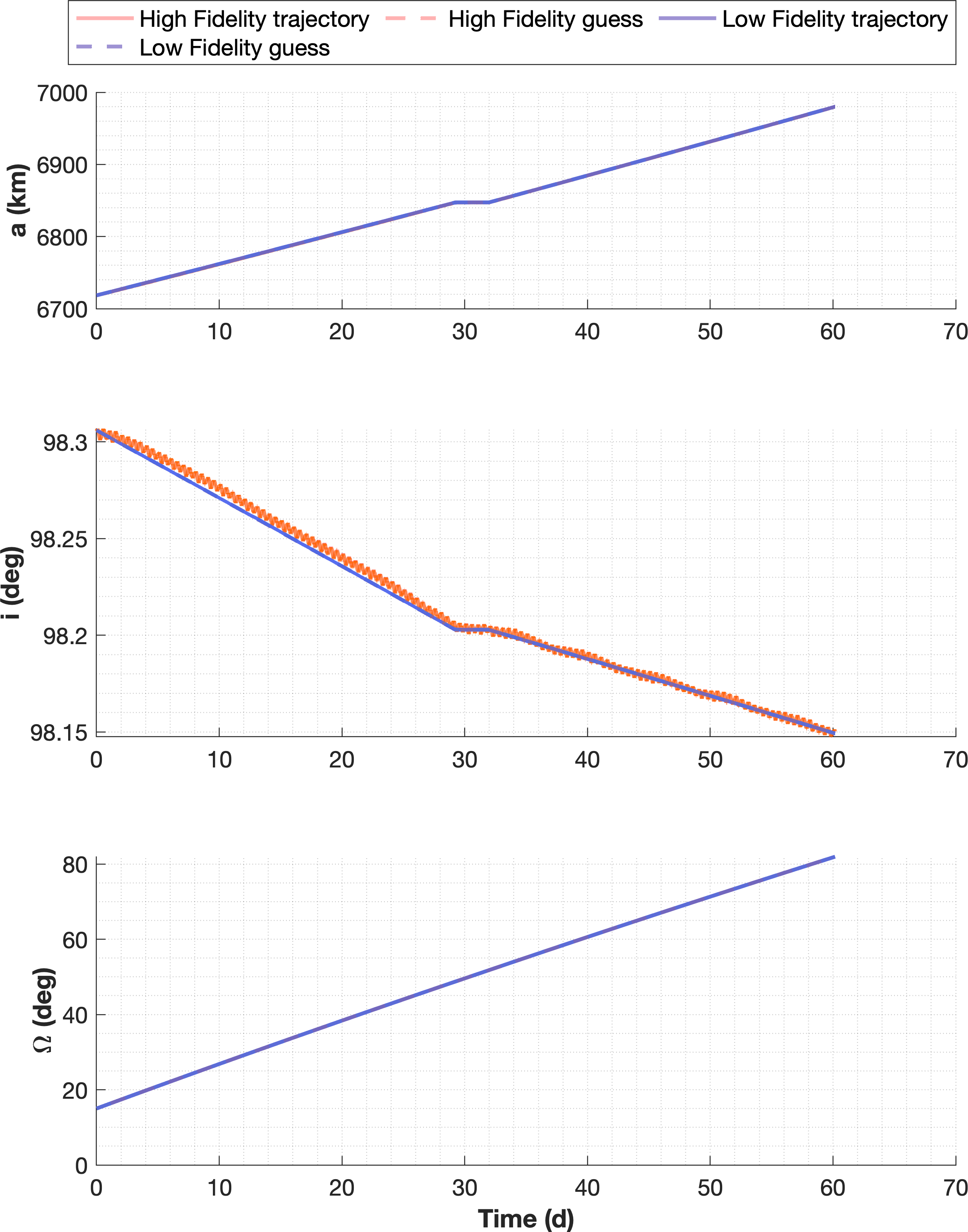}
    \caption{The fuel optimal solution with the high and low fidelity spacecraft propagations}
    \label{foptb11}
  \end{subfigure}
  \hfill
  \begin{subfigure}[b]{0.49\textwidth}
    \centering
    \includegraphics[width = \textwidth]{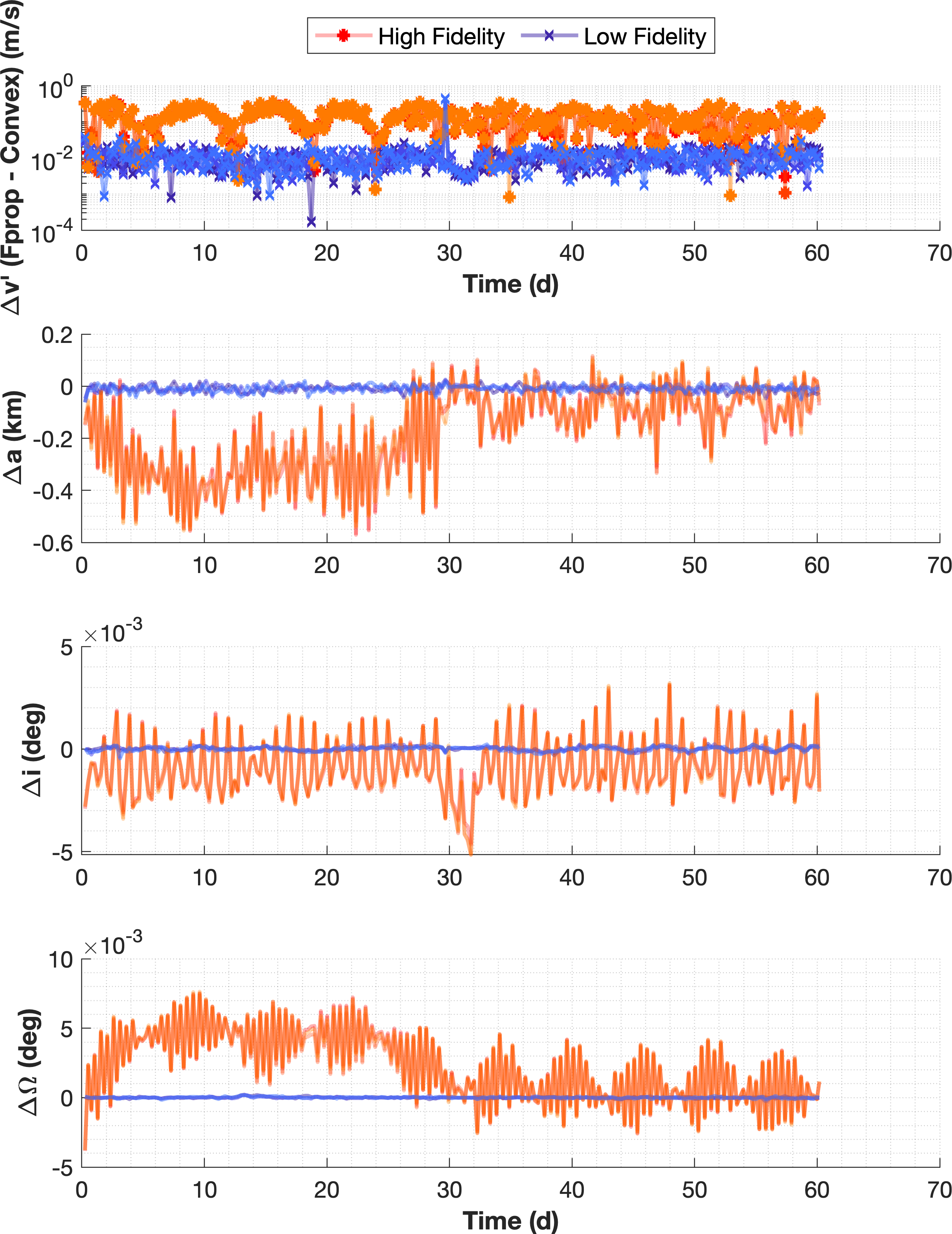}
   \caption{The evolution of $\Delta v'$, semi-major axis, inclination, and RAAN errors with respect to the reference trajectories over time}
    \label{foptb21}
  \end{subfigure}
  \caption{Results of the up leg fuel optimal implementation with high thrust errors }
  \label{fig:main3}
\end{figure}

The initial fuel optimal reference trajectory obtained from the PMDT is of  $\Delta v = \SI{148.9966}{\meter\per\second} $ and $TOF = \SI{60.1720}{\day}$. All simulations were done with $DC' = 0.4$, as this value was shown to give sufficient margin to mitigate significant deviations from the initial reference. Figures \ref{foptb1} and \ref{foptb2} shows the results with perfect thrust, while figures \ref{foptc11} and \ref{foptc21}
show the results with the first set of thrust errors. Figures \ref{foptb11} and \ref{foptb21} show the results with the second set of thrust errors. The numerical results are given in Table \ref{tR}. \par 
It can be seen that even with perfect thrust, the high-fidelity forward propagation results in a more significant fuel consumption than the low-fidelity one. This is expected, as more fuel will be required to account for deviations due to high-order perturbations, third-body effects, and solar radiation pressure. It was noted that the most significant deviations between high and low-fidelity propagation originated from the high-order gravitational perturbations, specifically the effect of zonal and sectorial harmonics on the inclination.
Simulations with thrust errors consume more fuel than ones without thrust errors for both propagations, which is also expected as propellant is required for error correction.  \par
\begin{figure}[hbt!]
  \centering
  \begin{subfigure}[b]{0.49\textwidth}
  \centering
    \includegraphics[width = \textwidth]{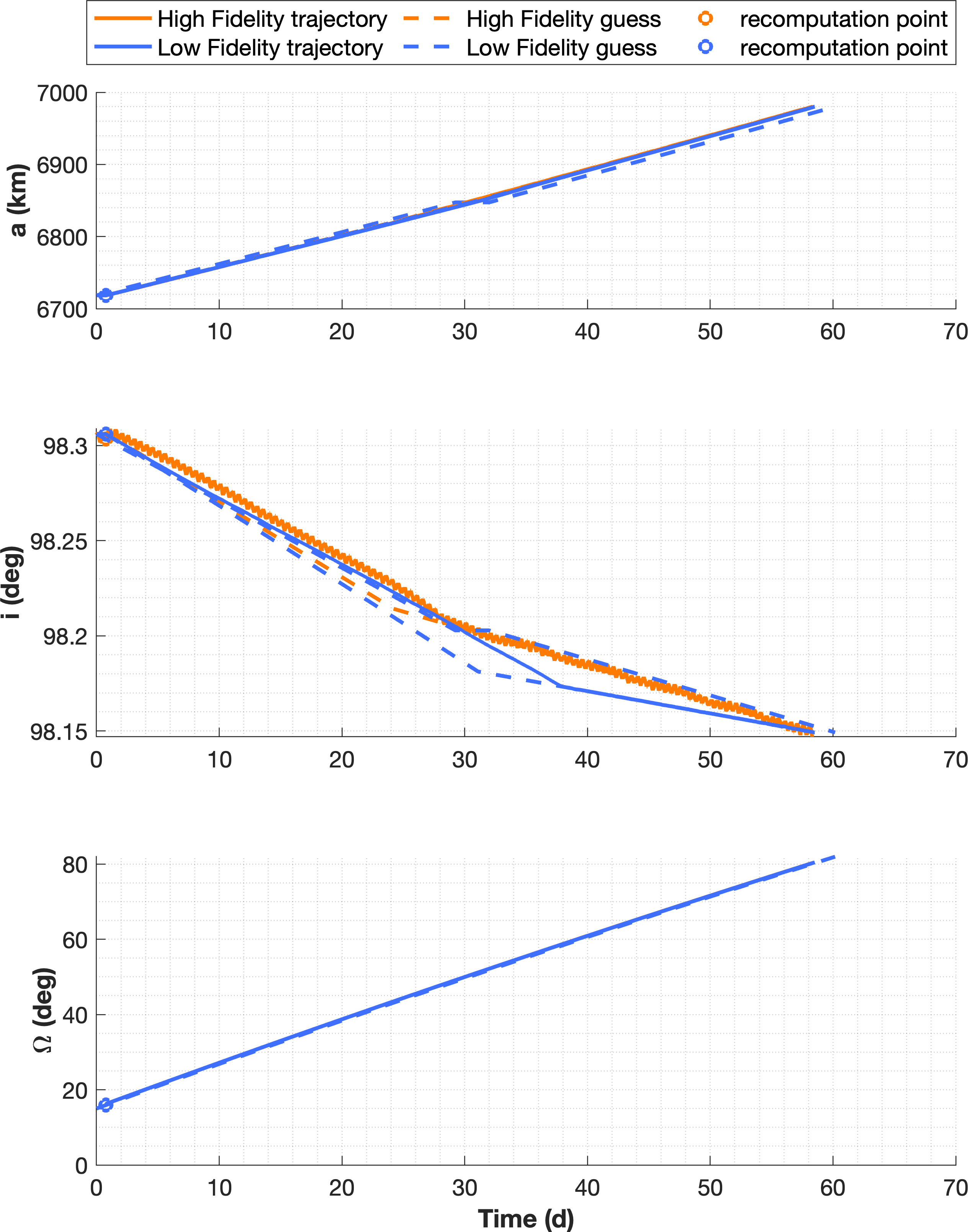}
    \caption{The fuel optimal solution with three low fidelity spacecraft propagations}
    \label{foptb12}
  \end{subfigure}
  \hfill
  \begin{subfigure}[b]{0.49\textwidth}
    \centering
    \includegraphics[width = \textwidth]{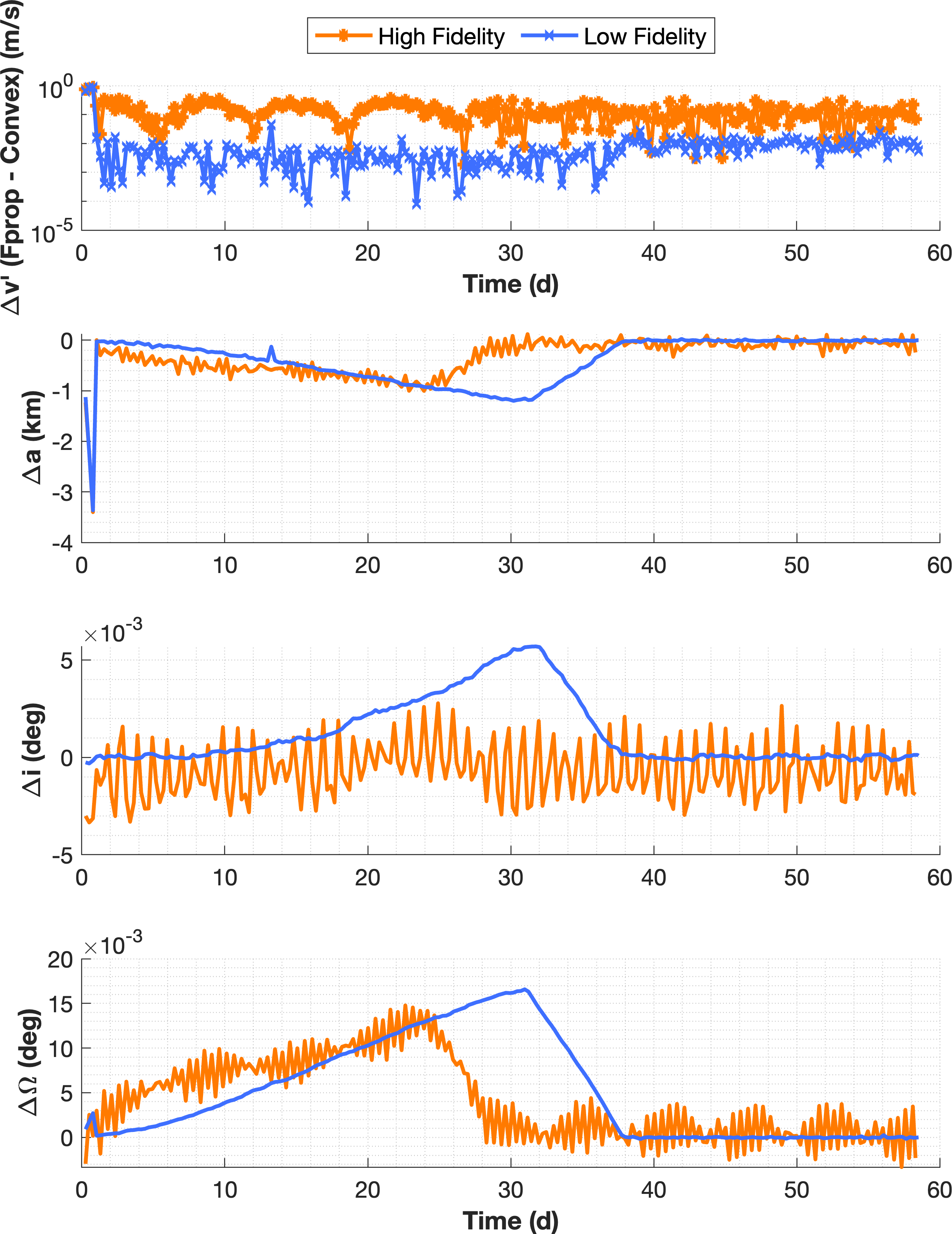}
    \caption{The evolution of $\Delta v'$, semi-major axis, inclination, and RAAN errors with respect to the reference trajectories over time}
    \label{foptb22}
  \end{subfigure}
  \caption{Results of the up leg fuel optimal implementation with high thrust errors and continuous misthrust for 15 orbits at the start}
  \label{fig:main4}
\end{figure}

The results shown in Figures \ref{foptb12} and \ref{foptb22} were generated with the high thrust errors for a case where the thrust is initially turned off for 15 orbits (3 MPC segments) to induce recomputations. Circles denote the recomputation points. It can be seen that the recomputations reduce the overall TOF for both propagations. While the fuel consumption of the low fidelity propagation is increased, that of the high fidelity propagation is slightly decreased. It should be noted that a recomputation may find trajectories that are more or less 'optimal' than the initial reference.

\begin{table}[hbt!]
\begin{tabular}{ccccccc}
Spacecraft Propagation & $\Delta a$ (km) & $\Delta i$ (deg) & $\Delta \Omega$ (deg) & TOF (d) & $\Delta v$ (m/s) & $\Delta v^{'}$ (m/s) \\ \hline 
\multicolumn{7}{c}{Perfect Thrust}                                                                                                      \\
Low Fidelity & -0.00040523 & -2.5314e-06 & 8.3604e-05 & 60.1721 & 151.997 & 0.02281 \\
High Fidelity & -0.050871 & 0.0012275 & -0.0019501 & 60.1721 & 161.1645 & 0.32583  \\ \hdashline
\multicolumn{7}{c}{Low thrust errors}                                                                                                      \\
 Low Fidelity (avg) & -0.0033709 & -9.8145e-07 & 9.7815e-05 & 60.1721 & 152.6132 & 0.028989 \\
 High Fidelity (avg) & -0.051819 & 0.0011801 & -0.0019695 & 60.1721 & 161.6146 & 0.32591 \\
\hdashline
\multicolumn{7}{c}{High thrust errors}                                                                                                      \\
Low Fidelity (avg) & -0.010719 & -9.3973e-06 & 9.176e-05 & 60.1721 & 153.1981 & 0.033133 \\
High Fidelity (avg) &-0.051597 & 0.0011507 & -0.0020179 & 60.1721 & 162.2609 & 0.32869\\ \hline 
\multicolumn{7}{c}{High thrust errors (with forced misthrust)}                                                                                                      \\
Low Fidelity (avg) &-0.0014997 & -1.8843e-06 & 8.3989e-05 & 58.5166 & 154.431 & 0.64575\\
High Fidelity (avg) &-0.24347 & -0.0023085 & -0.0019065 & 58.3046 & 160.9931 & 0.75872 \\ \hline 
\end{tabular}
\caption{Up leg fuel-optimal results}
\label{tR}
\end{table}

% \begin{table}[hbt!]
% \begin{tabular}{cccc}
%  Propagation & TOF (d) & $\Delta v$ (m/s) & $\Delta v^{'}$ (m/s) \\ \hline 
% \multicolumn{4}{c}{Perfect Thrust}                                                                                                      \\
% Low Fidelity  & 60.1721 & 151.997 & 0.02281 \\
% High Fidelity & 60.1721 & 161.1645 & 0.32583  \\ \hdashline
% \multicolumn{4}{c}{Low thrust errors}                                                                                                      \\
%  Low Fidelity (avg) & 60.1721 & 152.6132 & 0.028989 \\
%  High Fidelity (avg)  & 60.1721 & 161.6146 & 0.32591 \\
% \hdashline
% \multicolumn{4}{c}{High thrust errors}                                                                                                      \\
% Low Fidelity (avg) &  60.1721 & 153.1981 & 0.033133 \\
% High Fidelity (avg) & 60.1721 & 162.2609 & 0.32869\\ \hline 
% \multicolumn{4}{c}{High thrust errors (with forced misthrust)}                                                                                                      \\
% Low Fidelity (avg) & 58.5166 & 154.431 & 0.64575\\
% High Fidelity (avg) & 58.3046 & 160.9931 & 0.75872 \\ \hline 
% \end{tabular}
% \caption{Up leg fuel-optimal results}
% \label{tR}
% \end{table}

\subsection{Fuel optimal down leg from the target debris to the initial orbit}

Once the servicer spacecraft has reached the debris and rendezvoused with it, it is expected to bring the debris back to a 350 km altitude. At the start of the downward leg, the total mass is that of the dry mass of the servicer plus the remainder of the fuel mass and the debris mass. As such $m_0 = 796.2437 +2991 = \SI{ 3787.2}{\kilo\gram}$. Thirty days are allocated for the rendezvous procedure; hence the down leg start date is 00:00 UTC 20-Jun-2022. Table \ref{data2} shows the starting and target orbital parameters for the down leg. The drag coefficient and frontal area are again assumed to be 2.2 and $\SI{0.01}{\meter\squared}$, respectively. Note that during the downward leg, only the semi-major axis is tracked. The eccentricity is maintained to be quasi-circular by alternating the direction of the out-of-plane thrust as shown in Algorithm \ref{alg:3}. It is assumed that the Shepherd spacecraft will be able to reach and rendezvous with the debris once it is in a quasi-circular orbit at 350 km altitude, regardless of its other orbital parameters. 
% \begin{table} [h!]\centering\begin{tabular}{ccccc}Object & $a$ (km) & e & $i$ (deg) & $\Omega$ (deg) & $\omega$ (deg) & $\theta$ (deg) \\ \hline Initial (H-2A F-15) & 6987.0507 & 0.0042309 & 98.2219 & 108.8944 & 275.8823 & 64.4907 \\Target & 6728.1363 & 0.0042309 & 98.2219 & 108.8944 & 275.8823 & 64.4907 \\\hline \end{tabular} \caption{Orbital elements of initial debris and the target orbit at launch time of the down leg}  \label{data2} \end{table}

\begin{table} [h!]\centering\begin{tabular}{ccccccc}Object & $a$ (km) & e & $i$ (deg) & $\Omega$ (deg) & $\omega$ (deg) & $\theta$ (deg) \\ \hline Initial (H-2A F-15) & 6987.0507 & 0.0042309 & 98.2219 & 108.8944 & 275.8823 & 64.4907 \\Target & 6728.1363 & - & -& - & - & - \\\hline \end{tabular} \caption{Orbital elements of initial debris and the target orbit at launch time of the down leg}  \label{data2} \end{table}

\begin{figure}[hbt!]
  \centering
  \begin{subfigure}[b]{0.49\textwidth}
 \centering
    \includegraphics[width = \textwidth]{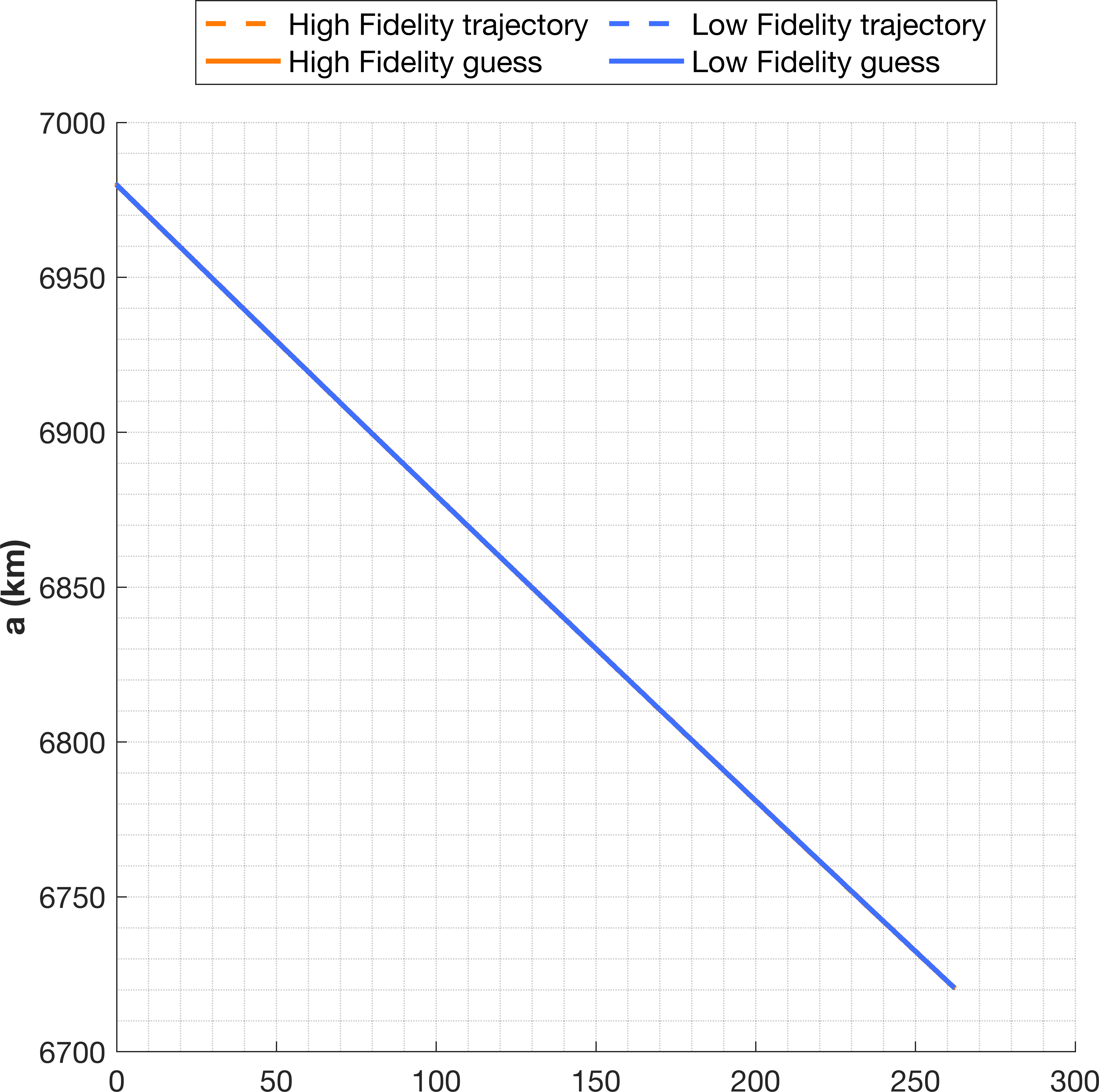}
    \caption{The fuel optimal solution with the high and low fidelity spacecraft propagations}
    \label{dlegne1}
  \end{subfigure}
  \hfill
  \begin{subfigure}[b]{0.49\textwidth}
   \centering
    \includegraphics[width = \textwidth]{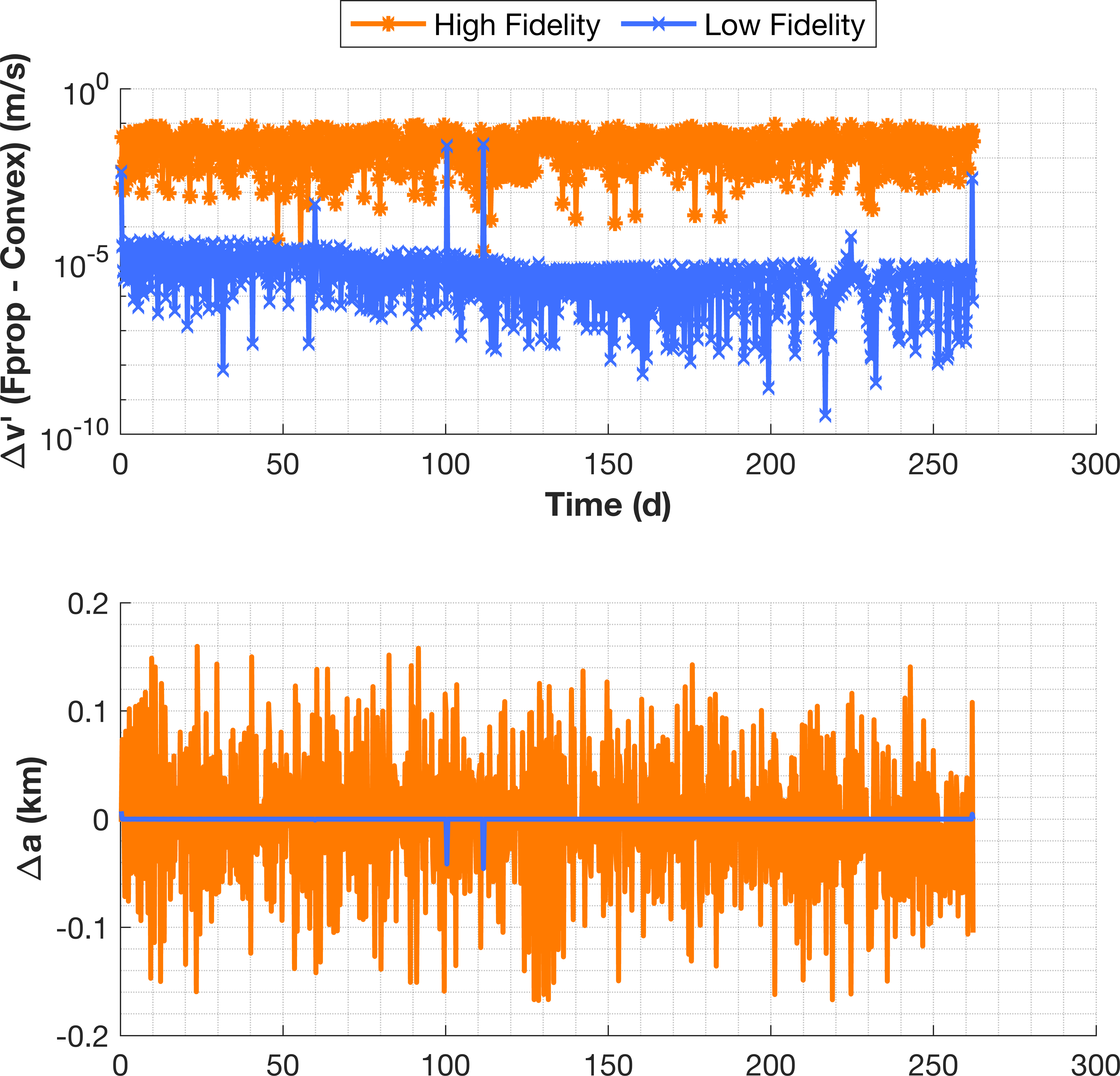}
    \caption{The evolution of $\Delta v'$ and  semi-major axis errors with respect to the reference trajectories over time}
    \label{dlegne2}
  \end{subfigure}
  \caption{Results of the down leg fuel optimal implementation without thrust errors}
  \label{fig:main4}
\end{figure}

\begin{figure}[hbt!]
  \centering
  \begin{subfigure}[b]{0.49\textwidth}
 \centering
    \includegraphics[width = \textwidth]{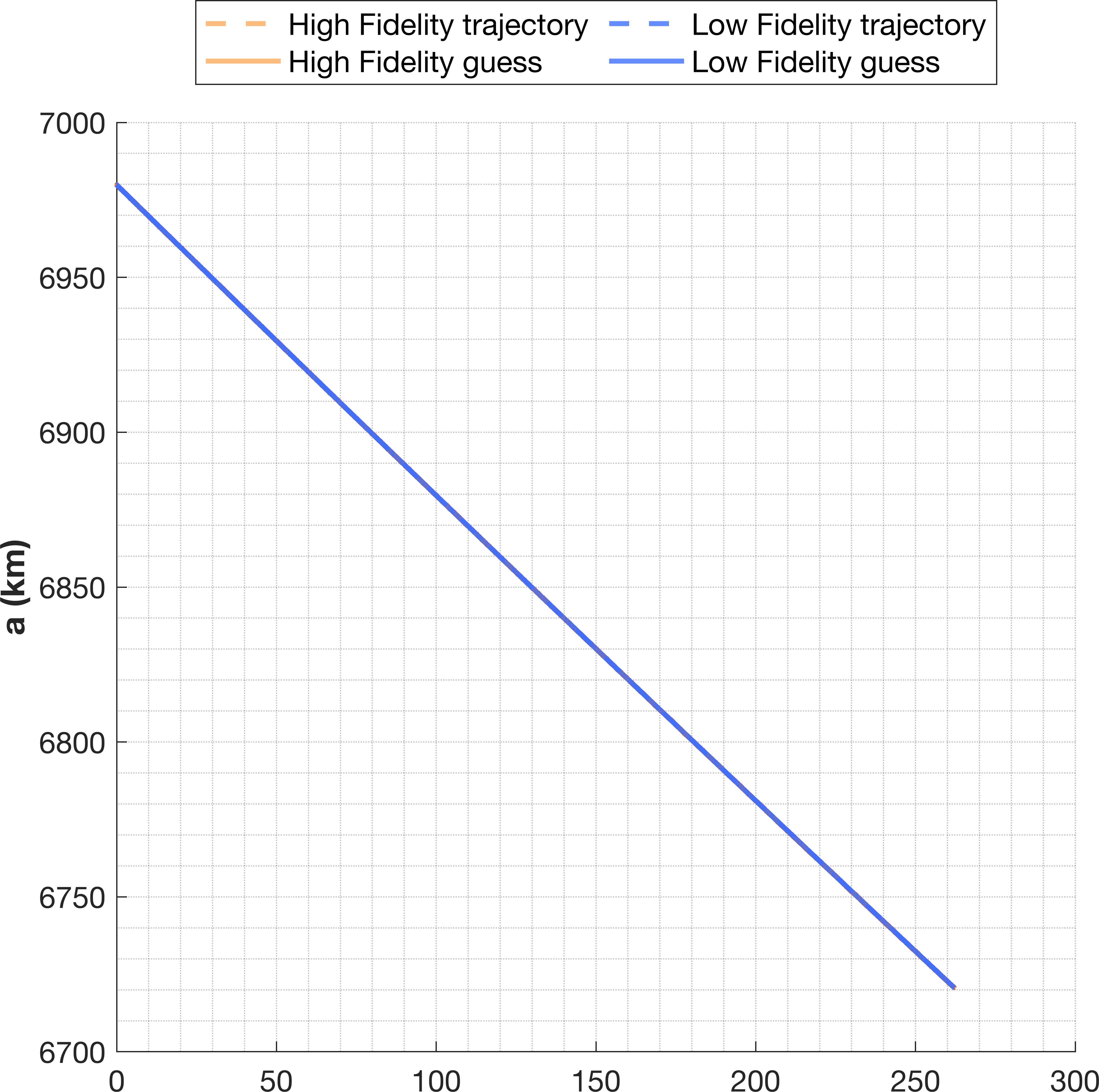}
    \caption{The fuel optimal solution with the high and low fidelity spacecraft propagations}
    \label{dleg1}
  \end{subfigure}
  \hfill
  \begin{subfigure}[b]{0.49\textwidth}
   \centering
    \includegraphics[width = \textwidth]{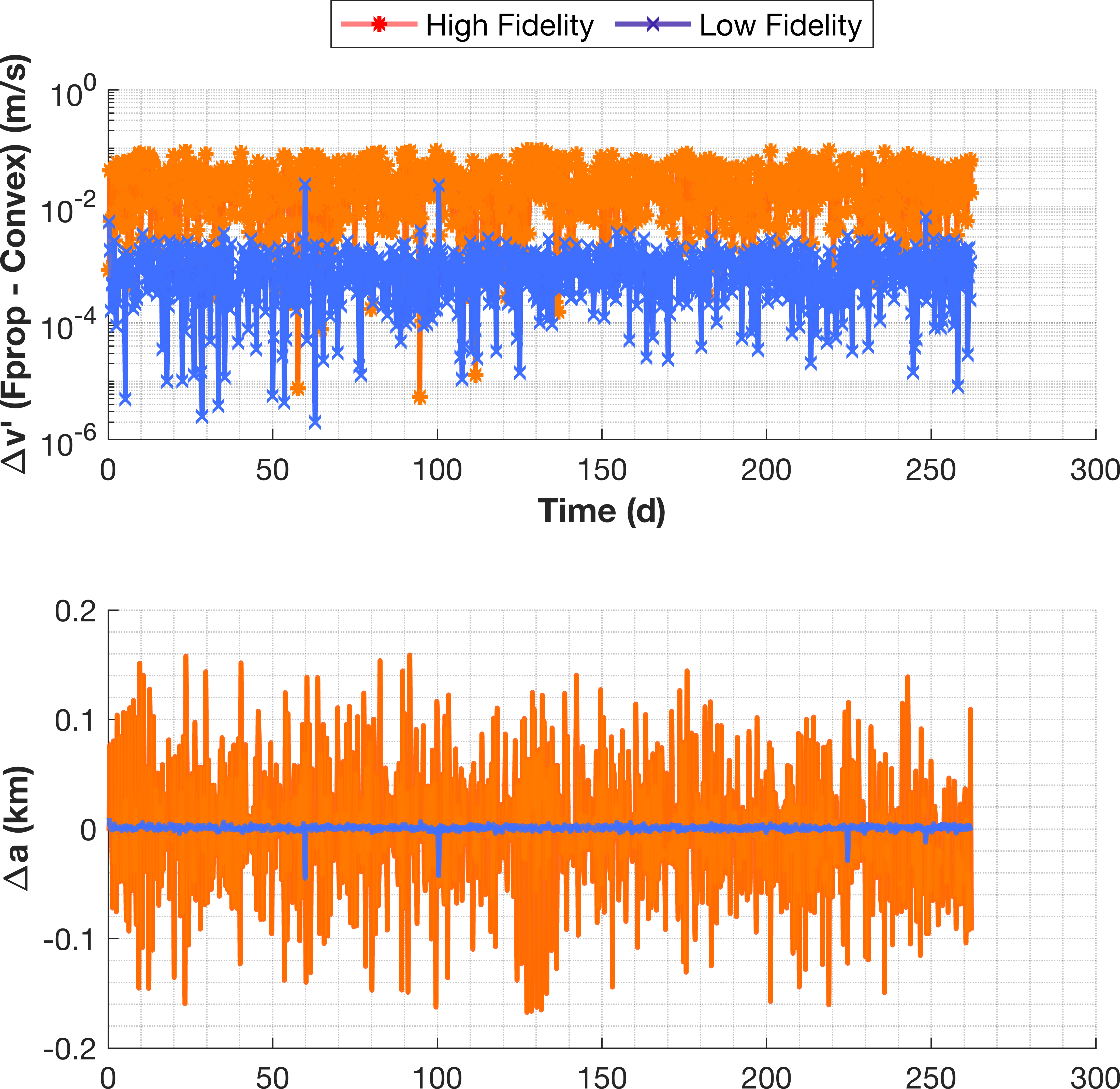}
    \caption{The evolution of $\Delta v'$ and semi-major axis over time}
    \label{dleg2}
  \end{subfigure}
  \caption{Results of the down leg fuel optimal implementation with low thrust errors}
  \label{fig:main5}
\end{figure}

\begin{figure}[hbt!]
  \centering
  \begin{subfigure}[b]{0.49\textwidth}
 \centering
    \includegraphics[width = \textwidth]{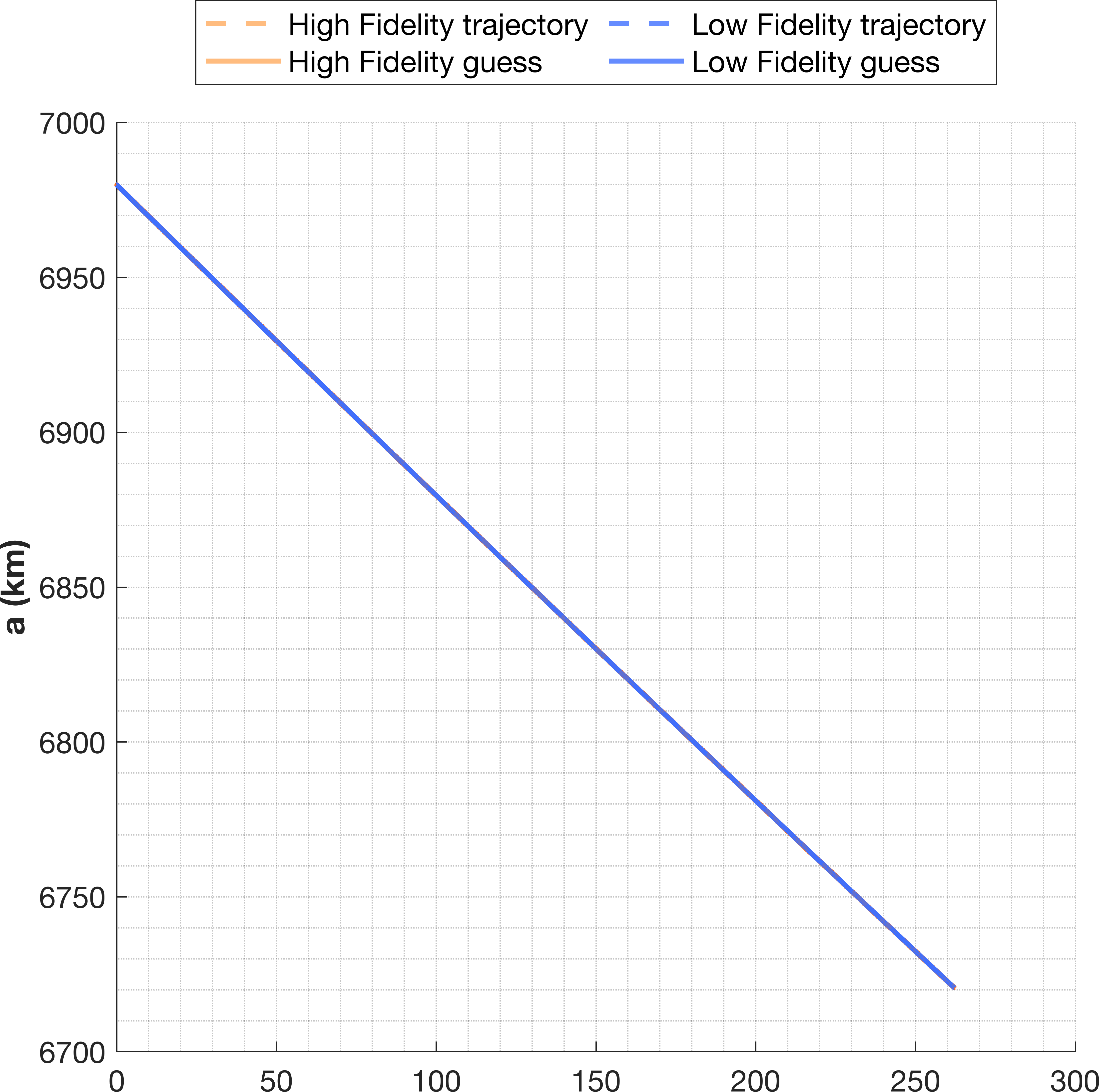}
    \caption{The fuel optimal solution with the high and low fidelity spacecraft propagations}
    \label{dleg1h}
  \end{subfigure}
  \hfill
  \begin{subfigure}[b]{0.49\textwidth}
   \centering
    \includegraphics[width = \textwidth]{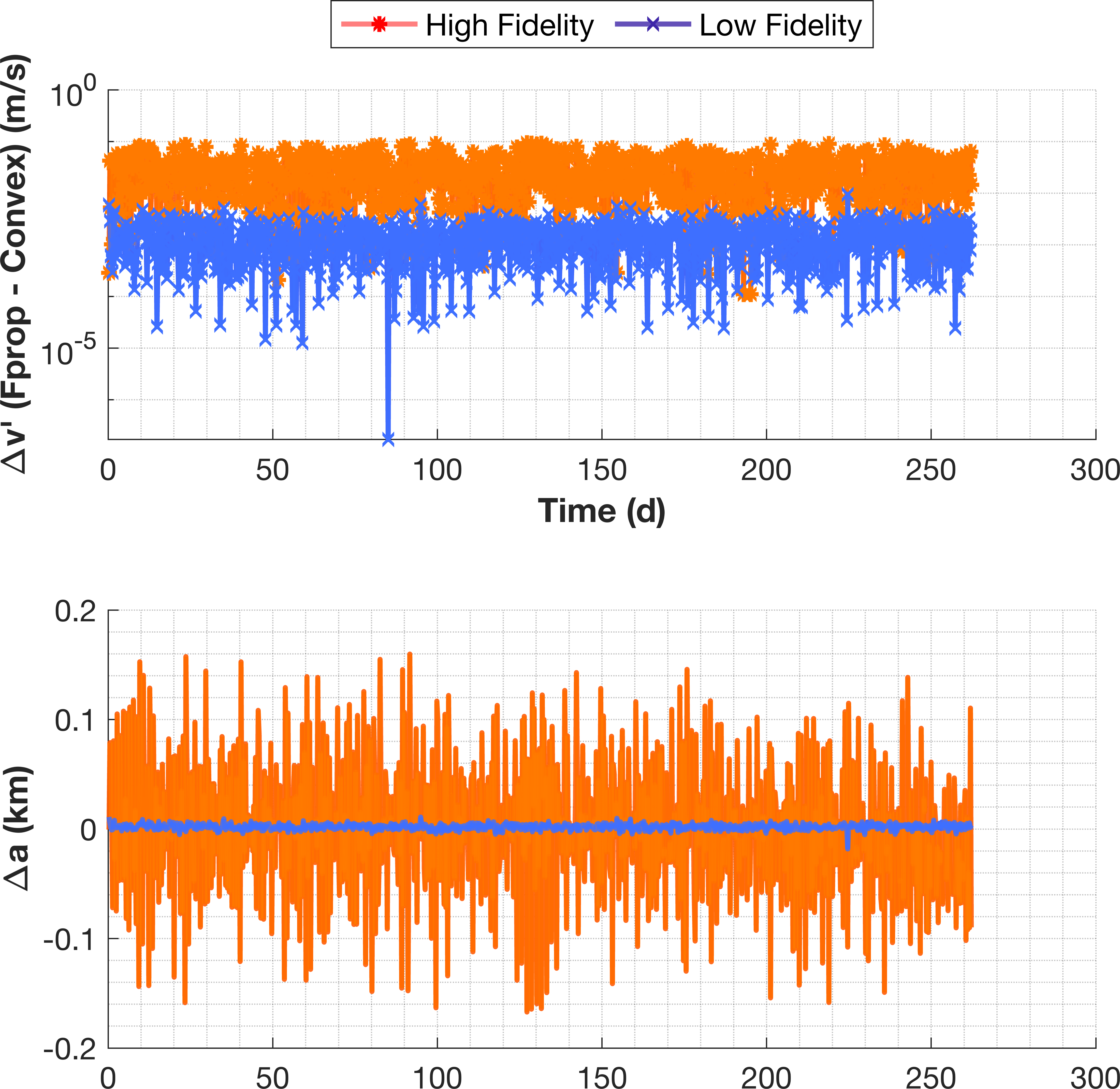}
    \caption{The evolution of $\Delta v'$ and semi-major axis over time}
    \label{dleg2h}
  \end{subfigure}
  \caption{Results of the down leg fuel optimal implementation with high thrust errors}
  \label{fig:main6}
\end{figure}

The initial fuel optimal reference trajectory obtained from the PMDT is of  $\Delta v = \SI{144.3437}{\meter\per\second} $ and $TOF = \SI{262.1590}{\day}$ for the down leg. All simulations were again done with $DC' = 0.4$. Figures \ref{dlegne1} and \ref{dlegne2} show the results with perfect thrust, while figures \ref{dleg1} and \ref{dleg2} show the results with the first set of thrust errors. Figures \ref{dleg1h} and \ref{dleg2h} show the results with the second set of thrust errors. The numerical results are given in Table \ref{td}. \par 
Unlike with the up leg, in the case of perfect thrust, the high-fidelity propagation requires less fuel than the low-fidelity one. This is likely because the PMDT uses the Harris-Priester (HP) atmospheric density model \cite{HATTEN2017571} while the high fidelity propagation uses the more accurate NRLMSISE-00 \cite{MORSELLI2014490}. NRLMSISE-00 estimates higher densities than the HP model at low altitudes \cite{compare}. Hence, the higher atmospheric drag in the high-fidelity propagation reduces the $\Delta v$ required to lower the spacecraft's altitude.  
It can be seen that the simulations with thrust errors consume more fuel than the case without thrust errors for both low and high-fidelity propagations, as expected. 

\begin{table}[hbt!]
\centering
\begin{tabular}{ccccc}
Spacecraft Propagation & $\Delta a$ (km) & TOF (d) & $\Delta v$ (m/s) & $\Delta v^{'}$ (m/s) \\ \hline 
\multicolumn{5}{c}{Perfect Thrust}                                                                                                      \\
Low Fidelity & -1.2481e-06 & 262.159 & 144.2285 & 0.0039075 \\
High Fidelity & -0.1049  & 262.159 & 140.8204 & 0.0013419 \\ \hdashline
\multicolumn{5}{c}{Low thrust errors}                                                                                                      \\
Low Fidelity (avg) & 0.0019581 & 262.159 & 144.7774 & 0.0054359 \\
High Fidelity (avg) & -0.092653 & 262.159 & 141.3518 & 0.00080687 \\
\hdashline
\multicolumn{5}{c}{High thrust errors}                                                                                                      \\
 Low Fidelity (avg) & 0.0032589  & 262.159 & 145.3054 & 0.0059401 \\
 High Fidelity (avg)& -0.089798 & 262.159 & 141.8686 & 0.00028516  \\
\hline 
\end{tabular}
\caption{Down leg fuel-optimal results}
\label{td}
\end{table}

\section{Conclusion}
This paper proposes a novel Convex Optimization-based Model Predictive Control (MPC) approach for providing guidance for Active Debris Removal (ADR) missions. The reference trajectory is updated in an MPC manner while tracking is conducted using convex optimization. For each optimization segment, boundary conditions are only imposed at the start, and the end, so the trajectory can follow a pre-optimized path closely while avoiding unnecessary fuel expenditure by mitigating the need to match the optimized trajectory exactly. \par 
This work stems from the work done in [\citenum{ADRMW}], where preliminary optimal trajectories were developed for debris removal. The reference trajectories for the guidance to follow are obtained using the same methodology. \par 
The key contribution of this work lies in eliminating successive convexification of nonconvex dynamics, which has been a source of inaccuracies in the past. By employing near-linear Generalized Equinoctial Orbital Elements and the highly accurate PMDT as an initial guess and a reference, the convex guidance achieves accurate solutions without requiring iterative computations. The MPC recomputes an optimal trajectory upon encountering large deviations from the optimal path.  A single parameter-$\Delta v'$-is developed to gauge the accuracy of the trajectory relative to the reference to determine if a recomputation is warranted. Furthermore, the reference is calculated at a lower duty cycle than the real duty cycle to allow for a thrust margin that could be used to counteract minor trajectory deviations, prolonging the validity of a reference. \par 
Results are provided for a mission that starts from a 350 km parking orbit to the tonne class rocket body H-2A (F15), which is then collected and brought back to a 350 km orbit. The spacecraft is propagated with high and low-fidelity dynamics and thrust errors to emulate reality.  The results show that the spacecraft can closely follow the optimized reference despite thrust errors, even under high-fidelity dynamics. The spacecraft can also recompute an optimal trajectory if it deviates significantly from the original reference.\par 

\clearpage
\section{Acknowledgment}
This work was partially supported by the Ministry of Business, Innovation, and Employment (MBIE) study: Astroscale/ Rocket Lab/ Te Pūnaha Ātea-Space Institute Active Debris Removal Study.
%\section*{Acknowledgements}

\bibliographystyle{AAS_publication}   % Number the references.

\bibliography{library.bib}   % Use references.bib to resolve the labels.

\end{document}